\documentclass[reqno,12pt]{amsart}
\usepackage{latexsym}
\usepackage{amsmath,amssymb}
\usepackage{amsthm}
\usepackage{amsfonts}
\numberwithin{equation}{section}
\newtheorem{definition}{Definition}[section]
\newtheorem{theorem}{Theorem}[section]
\newtheorem{lemma}{Lemma}[section]

\newtheorem{proposition}{Proposition}[section]
\newtheorem{remark}{Remark} [section]

\date{}
\title{On One-dimensional  Compressible Navier-Stokes Equations with Degenerate Viscosity and
Constant State at Far fields}

\author[Changsheng Dou, Quansen Jiu]{}

\begin{document}
\maketitle

%-----------------------------------------------------------------------------------

   \centerline{\scshape Changsheng Dou }

\medskip
{\footnotesize
  \centerline{Institute of Applied Physics and Computational Mathematics}
  \centerline{Beijing  100088, P. R. China}
   \centerline{  \it Email: douchangsheng@163.com}
   }
\vspace{3mm}

 \centerline{\scshape Quansen Jiu$^{\star}$ }

\medskip
{\footnotesize
  \centerline{School of Mathematical Sciences, Capital Normal University}
  \centerline{Beijing  100048, P. R. China}
   \centerline{  \it Email: jiuqs@mail.cnu.edu.cn}
   }
\footnotetext{

$^\star$ Partially supported by NSFC grants ( No. 11171229).}

\pagestyle{myheadings}
\thispagestyle{plain}\markboth{\small{Changsheng Dou, Quansen Jiu }}
{\small{One-dimensional Compressible NS Equations }}

\vspace{3mm}
\begin{center}
\textbf{Abstract}
\end{center}
In this paper, we are concerned with
 the Cauchy problem for one-dimensional compressible isentropic Navier-Stokes equations with density-dependent
viscosity $\mu(\rho)=\rho^\alpha (\alpha>0)$ and  pressure
$P(\rho)=\rho^{\gamma}\ (\gamma>1)$.  We will establish the global
existence and asymptotic behavior of  weak solutions for  any
$\alpha>0$ and $\gamma>1$ under the assumption that the density
function keeps a constant state at far fields. This enlarges the
ranges of $\alpha$ and $\gamma$ and improves the previous results
presented by Jiu and Xin.  As a
 result, in the case that  $0<\alpha<\frac12$,
 we obtain the large time  behavior of the
strong solution  obtained  by  Mellet and Vasseur when the solution
has a lower bound (no vacuum).

\section{Introduction}

Consider the one-dimensional compressible Navier-Stokes equations
with density-dependent viscosity coefficients
\begin{equation}\label{2.2.2}
\left\{\begin{array}{llll}
 \rho_t+(\rho u)_x = 0 ,\\
 (\rho u)_t+(\rho u^2 +P(\rho))_x=(\mu(\rho) u_x)_x.\\
 \end{array}\right.
\end{equation}
Here, $\rho(x, t)$ and $u(x, t)$ stand for the fluid density and
velocity respectively. For simplicity, the pressure term $P(\rho)$
and the viscosity coefficient $\mu(\rho)$ are assumed to be
\begin{equation}\label{2.2.3}
P(\rho)=\rho^{\gamma}\ (\gamma>1),\ \ \ \ \mu(\rho)=\rho^{\alpha}.
\end{equation}
The initial data is imposed as
\begin{equation}\label{2.2.4}
(\rho, \rho u)|_{t=0}=(\rho_0, m_0).
\end{equation}

When the viscosity $\mu(\rho)$ is a positive constant, there has
been a lot of investigations on the compressible Navier-Stokes
equations, for smooth initial data or discontinuous initial data,
one-dimensional or multidimensional problems (see
\cite{KS1977,Serre1986,Hoff1987,Nash1962,MN1980,Danchin2000,Lions1998,FNP2001},
 and the references
therein). However, the studies in Hoff $\&$ Serre \cite{HS1991}, Xin
\cite{Xin1998}, Liu, Xin $\&$ Yang \cite{LXY1998} show that the
compressible Navier-Stokes equations with constant viscosity
coefficients behave singularly in the presence of vacuum. In
\cite{LXY1998}, Liu, Xin and Yang introduced the modified
compressible Navier-Stokes equations with density-dependent
viscosity coefficients for isentropic fluids. Actually, in deriving
the compressible Navier-Stokes equations from the Boltzmann
equations by the Chapman-Enskog expansions, the viscosity depends on
the temperature, and correspondingly depends on the density for
isentropic cases. Meanwhile, an one-dimensional viscous Saint-Venant
system for shallow water, which was derived rigorously by
Gerbeau-Perthame \cite{GP2001}, is expressed exactly as
(\ref{2.2.2}) with $\mu(\rho)=\rho$  and $\gamma=2$.

When the viscosity $\mu(\rho)$ depends on the density, there are a
large number of literatures on mathematical studies on
(\ref{2.2.2}). One-dimensional case is referred to
\cite{LXY1998,Jiang1998,OMM2002,JXZ2005,YYZ2001,YZ2002,VYZ2003,LLX2008,JX2008,DJ2010}
and references therein. In  \cite{JX2008} Jiu and Xin studied the
global existence and large time behavior of weak solutions of the
Cauchy problem to (\ref{2.2.2}) with $\mu({\rho})={\rho}^{\alpha}$
($\alpha>1/2$) under some restrictions of $\alpha$ and $\gamma$. The
vacuum or non-vacuum constant states at far fields are permitted in
\cite{JX2008}. Recently, based on \cite{JX2008}, Guo etc
\cite{WHG2011} studied the global existence and large time behavior
of weak solutions  of the Cauchy problem to (\ref{2.2.2}) under the
assumptions of $0<\alpha<\frac{1}{2}$ and ${\rho}_0\in L^1(R)$. If
the far fields hold different ends, the asymptotic stability of
rarefaction waves for the compressible isentropic Navier-Stokes
equations (\ref{2.2.2}) with $\mu({\rho})={\rho}^{\alpha}
(\alpha>\frac12)$ was studied by Jiu, Wang and Xin in \cite{JWX2011}
in which the rarefaction wave itself has no vacuum, and in
\cite{JWX2011-2} in which the rarefaction wave connects with the
vacuum. In  \cite{MV2008},  Mellet and Vasseur showed that if
$0<\alpha<\frac{1}{2}$ and the initial datum are regular with a
positive lower bound (no vacuum), there exists a global and unique
strong solution of the Cauchy problem to (\ref{2.2.2}). However, the
a priori estimates obtained in \cite{MV2008} depends on the time
interval and hence does not yield the time-asymptotic behavior of the
solutions.

In this paper, we will study the global existence and asymptotic
behavior of weak solutions for  any $\alpha>0$ and $\gamma>1$ under
the assumption that the density function keeps a constant state at
far fields.  We will apply the similar approaches as in
\cite{JWX2011} to obtain an uniform (in time) entropy estimate (see
Section 3). This type of entropy estimate was observed first in
\cite{Kanel1968} for the one-dimensional case and later established
in \cite{BD2003,BD2004,BDL2003} for multi-dimensional cases.   The
key points in our proof are to obtain the uniform upper bound of the
density and to obtain the lower bound of the density of the
approximate solutions by using the uniform entropy estimate. To do
that, different ranges of $\alpha$ and $\gamma$ will be discussed
respectively and the elaborate estimates  will be given.  Our
results  relax the restrictions of $\alpha$ and $\gamma$ presented
in \cite{JX2008}.  In the case that  $0<\alpha<\frac12$,
 we obtain the large time asymptotic behavior of the
strong solution  obtained  by  Mellet and Vasseur when the solution
has a lower bound (no vacuum). Moreover, in the case that
$\alpha>\frac 12$, the vacuum is permitted and we study the
existence and large time behavior in the framework of weak
solutions.

The organization of the paper  is  as follows. In Section 2, we
state some preliminaries and  main results. In Section 3, we give
proofs of uniform entropy estimates. Based on these, lower and upper
bounds of the density  to the approximate solutions will be  shown.
In Section 4, we give a sketch of  proof of main results.

\section{Preliminaries and  Main results }

We first give the assumptions of the initial data and the definition
of weak solutions.

Define the  pressure potential $\Psi(\rho,\bar{\rho})$ as
\begin{equation}\label{2.2.6'}
\Psi(\rho, \bar{\rho})=\int_{\bar{\rho}}^{\rho}
\frac{p(s)-p(\bar{\rho})}{s^2} ds
=\frac{1}{(\gamma-1)\rho}[{\rho}^\gamma -{\bar{\rho}}^\gamma -\gamma
{\bar{\rho}}^{\gamma-1}(\rho-\bar{\rho})].
\end{equation}

We assume that there exists a constant $\bar{\rho}\geq 0$ such that
\begin{equation}\label{2.2.6}
\rho_0\Psi(\rho_0, \bar{\rho})\in L^1(R) \ \ {\rm if } \ \bar\rho>
0,\ \ {\rm and}\ \ \rho_0\in L^1(R) \ \ {\rm if}\ \ \bar\rho=0.
\end{equation}

Moreover, we assume that the initial data  satisfy
\begin{equation}\label{2.2.5}
\left\{\begin{array}{llll} {\rho}_0(x)> 0\ \ {\rm if}\ 0<\alpha\le
\frac12,\ \
\ \ {\rho}_0(x)\ge 0\ \ {\rm if}\ \alpha> \frac12;\\[3mm]
({{\rho}_0}^{\alpha-1/2})_x \in L^2(R)\ \
{\rm if }\ \alpha>0 \ \ {\rm and}\ \alpha\neq\frac12,\ \\[3mm] (\log
{{\rho}_0})_x\in L^2(R)\
\ {\rm if}\ \alpha=\frac{1}{2};\\[3mm]
 \frac{|m_0|^2}{\rho_0}\in L^1(R),\ \
\frac{|m_0|^{2+\delta}}{{\rho_0}^{1+\delta}}\in L^1(R),
\end{array}\right.
\end{equation}
where  $0<\delta<1$ is any fixed number which may be small.

\begin{remark} By assumptions \eqref{2.2.6}-\eqref{2.2.5}, the initial data $\rho_0(x)$  is actually continuous and bounded.
And in the case that $0<\alpha\le 1/2$, the restriction
${\rho}_0(x)> 0, x\in R$
 can be derived from other conditions of
\eqref{2.2.6}-\eqref{2.2.5}. However, in this case, the initial
density still can appear vacuum at infinity, i.e.,
$\lim_{|x|\to\infty} \rho_0(x)=0$.
\end{remark}

The  weak solutions to $(\ref{2.2.2})-(\ref{2.2.4})$ with the far
fields $\bar\rho\ge 0$ are defined as:

\begin{definition}\label{def2.2.1} For any $T>0$, a pair $(\rho, u)$ is said to be a weak solution to
$(\ref{2.2.2})-(\ref{2.2.4})$  if\\
(1) $\rho\geq 0$ a.e., and
\begin{eqnarray}
&&\rho-\bar\rho\in  L^\infty(0,T; L^\infty(R)),\\[3mm]
&& \rho\Psi(\rho, \bar{\rho})\in L^1(R), \sqrt{\rho} u \in L^{\infty}(0, T; L^2(R)),\nonumber\\[3mm]
&& (\rho^{\alpha-\frac12})_x\in L^\infty(0,T; L^2(R))\ \ {\rm if}\
 \ 0<\alpha\neq\frac12,\  \ \nonumber \\[3mm]
 &&(\log{\rho})_x \in L^{\infty}(0, T; L^2(R))\ \  {\rm if}\ \
 \alpha=\frac{1}{2};
\end{eqnarray}
(2) For any $t_2\geq t_1\geq 0$ and $\xi\in C^1(R\times[t_1, t_2])$, the mass equation $(\ref{2.2.2})_1$ holds in the following sense.\\
\begin{equation}\label{2.2.7}
{\int_{R}\rho\zeta
dx}|_{t_1}^{t_2}={\int_{{t_1}}^{t_2}}{\int_{R}}(\rho{\zeta}_t+\rho u
{\zeta}_x)dxdt;
\end{equation}
(3) For any $\psi\in C_0^{\infty}(R\times[0, T))$, the momentum equation holds in the following sense.\\
\begin{equation}\label{2.2.8}
{\int_{R}}m_0 \psi(0, \cdot)
dx+{\int_{t_1}^{t_2}}{\int_{R}}[\sqrt{\rho}(\sqrt{\rho}u)\psi_t+((\sqrt{\rho}u)^2+\rho^{\gamma})\psi_x]dxdt
+<\rho^{\alpha}u_x, \psi_x>=0;
\end{equation}
where the diffusion term make sense in the following equalities:\\
when $0<\alpha\neq\frac12$,
\begin{equation}\label{2.2.9'}
<\rho^{\alpha}u_x,
\psi_x>=-\int_{0}^{T}\int_{R}\rho^{\alpha-\frac{1}{2}}\sqrt{\rho}u\psi_x
dxdt-\frac{2\alpha}{2\alpha-1}\int_0^T\int_R
({\rho}^{\alpha-\frac{1}{2}})_x\sqrt{\rho}u\psi dxdt,
\end{equation}
when $\alpha=\frac{1}{2}$,
\begin{equation}\label{2.2.9}
<\rho^{\frac{1}{2}}u_x,
\psi_x>=-\int_{0}^{T}\int_{R}\sqrt{\rho}u\psi_x
dxdt-\frac{1}{2}\int_0^T\int_R (\log{\rho})_x\sqrt{\rho}u\psi dxdt.
\end{equation}
\end{definition}

Before we state our main results, we review the existence results
obtained in \cite{JX2008} as follows.

\begin{proposition}\label{prop1.1} (\cite{JX2008}, $\bar\rho=0$) Let $\gamma>1$ and $\alpha>\frac
12$.  Suppose that $(\ref{2.2.6})$ and $(\ref{2.2.5})$ hold. If
$\bar\rho=0$, then the Cauchy problem $(\ref{2.2.2})-(\ref{2.2.4})$
admits a global weak solution $(\rho(x,t), u(x,t))$ satisfying
\begin{eqnarray}
\rho\in  C(R\times (0,T)).\label{th11-1}
\end{eqnarray}
Moreover, one has
\begin{eqnarray}
&&\sup_{t\in[0,T]}\int_R \rho dx+\max_{(x,t)\in R\times[0,T]}\rho\le
C, \label{th11-2}\\ &&\sup_{t\in[0,T]}\int_R
(|\sqrt{\rho}u|^2+(\rho^{\alpha-\frac12})_x^2+\frac{1}{\gamma-1}\rho^\gamma
dx+\int_0^T\int_R [(\rho^{\frac{\gamma+\alpha-1}{2}})_x]^2 dxdt\le
C,\label{th11-3}
\end{eqnarray}
where $C$ is an absolute constant.
\end{proposition}

\begin{proposition}\label{prop1.2} (\cite{JX2008},$\bar\rho>0$)
Let $\alpha$ and $\gamma$ satisfy
\begin{equation}\label{prop1.2-2}
\gamma>1, \frac12<\alpha\le\frac32\ \ {\rm or}\ \  \gamma\ge
2\alpha-1, \alpha > \frac32.
\end{equation}
 Suppose that $(\ref{2.2.6})$ and $(\ref{2.2.5})$
hold. If $\bar\rho>0$, then the Cauchy problem
$(\ref{2.2.2})-(\ref{2.2.4})$ admits a global weak solution
$(\rho(x,t), u(x,t))$ satisfying
\begin{eqnarray}
\rho\in  C(R\times (0,T)).\label{th12-1}
\end{eqnarray}
Moreover, one has
\begin{eqnarray}
&&\sup_{t\in[0,T]}\int_R |\rho-\bar\rho|^2 dx+\max_{(x,t)\in
R\times[0,T]}\rho\le C, \label{th12-2}\\
&&\sup_{t\in[0,T]}\int_R
(|\sqrt{\rho}u|^2+(\rho^{\alpha-\frac12})_x^2+\frac{1}{\gamma-1}(\rho^\gamma-(\bar\rho)^\gamma-\gamma(\bar\rho)^{\gamma-1}(\rho-\bar\rho)))
dx\label{prop1.2-1}\\
&&+\int_0^T\int_R [(\rho^{\frac{\gamma+\alpha-1}{2}})_x]^2 dxdt\le
C,\label{th12-3}
\end{eqnarray}
where $C$ is an absolute constant.
\end{proposition}

\begin{remark}
It should be noted that in Proposition \ref{prop1.2}, the
restrictions of $\gamma$ and $\alpha$ \eqref{prop1.2-2} are
different from ones presented in Theorem 2.2 in \cite{JX2008}. This
is due to  that in \cite{JX2008}, instead of
$\rho-\bar\rho\in L^\infty(0,\infty;L^1(R))$, one should use  the fact that
$\rho-\bar\rho\in L^\infty(0,\infty;L^2(R))$ which follows from the
estimate of \eqref{prop1.2-1}.
\end{remark}

Our main results are  as follows.
\begin{theorem}\label{thm2.2.1} Let $\gamma>1,\ 0<\alpha\le \frac{1}{2}$ and assume that $(\ref{2.2.6})-(\ref{2.2.5})$ hold. Then for any $T>0$,
 the Cauchy problem $(\ref{2.2.2})-(\ref{2.2.4})$ admits a global weak solution $(\rho(x, t), u(x, t))$ in $R\times
 (0,T)$ satisfying

(1)
\begin{equation}\label{2.2.10}
 \rho\in C(R\times (0,T)), \ \ \rho(x,t)\ge 0, \ \ (x,t)\in R\times
(0,T);
\end{equation}

(2)
\begin{eqnarray}
 &&\sup_{t\in[0,T]}\int_{R}|\rho-\bar{\rho}|^2dx+\max_{(x,t)\in
R\times [0,T]}\rho\leq C,\ \ \ {\rm if}\ \ \bar{\rho}>0, \label{2.2.11}\\[3mm]
&&\sup_{t\in[0,T]}\int_{R}\rho dx+\max_{(x,t)\in R\times
[0,T]}\rho\leq C,\ \ \ {\rm if}\ \  \bar{\rho}=0;
\end{eqnarray}

(3) When $0<\alpha<\frac12$, one has
\begin{equation}\label{2.2.12}
\begin{aligned}
&\sup_{t\in[0,T]}\int_{R}(|\sqrt{\rho}u|^2+{(\rho^{\alpha-1/2})_x}^2+\frac{1}{\gamma-1}(\rho^\gamma-(\bar{\rho})^\gamma
-\gamma(\bar{\rho})^{\gamma-1}(\rho-\bar{\rho})))dx\\[3mm]
&+\int_{0}^{T}\int_{R}([(\rho^{\frac{\gamma+\alpha-1}{2}})_x]^2+{\Lambda(x,t)}^2)dxdt\leq
C,
\end{aligned}
\end{equation}
where $C$ is an absolute constant which only depends on the initial
data, and $\Lambda(x,t)\in L^2(R\times(0,T))$ is a function which
satisfies
$$
\int_{0}^{T}\int_{R}\Lambda\psi dxdt
=-\int_{0}^{T}\int_{R}{\rho}^{\alpha-1/2}\sqrt{\rho}u\psi_x
dxdt-\frac{2\alpha}{2\alpha-1}\int_0^T\int_R
{\rho}^{\alpha-1/2}_x\sqrt{\rho}u\psi dxdt;
$$

(4) When $\alpha=\frac12$, one has
\begin{equation}
\begin{aligned}\label{2.2.16}
&\sup_{t\in[0,T]}\int_{R}(|\sqrt{\rho}u|^2+{(\log{\rho})_x}^2+\frac{1}{\gamma-1}(\rho^\gamma-(\bar{\rho})^\gamma
-\gamma(\bar{\rho})^{\gamma-1}(\rho-\bar{\rho})))dx\\
&+\int_{0}^{T}\int_{R}([(\rho^{\frac{\gamma-\frac{1}{2}}{2}})_x]^2+{\Lambda(x,t)}^2)dxdt\leq
C,
\end{aligned}
\end{equation}
where $C$ is an absolute constant which just depends on the initial
data, and $\Lambda(x,t)\in L^2(R\times(0,T))$ is a function which
satisfies
$$
\int_{0}^{T}\int_{R}\Lambda\psi dxdt
=-\int_{0}^{T}\int_{R}\sqrt{\rho}u\psi_x
dxdt-\frac{1}{2}\int_0^T\int_R (\log{\rho})_x\sqrt{\rho}u\psi dxdt.
$$
\end{theorem}
\begin{theorem}\label{thm2.2.1'} Let $\alpha$ and $\gamma$ satisfy
\begin{equation}\label{ag}
\alpha>\frac12,\ \gamma>1.
\end{equation}
Suppose that $(\ref{2.2.6})-(\ref{2.2.5})$ hold. If $\bar{\rho}>0$,
then for any $T>0$,
 the Cauchy problem $(\ref{2.2.2})-(\ref{2.2.4})$ admits a global weak solution $(\rho(x, t), u(x, t))$ in $R\times
 (0,T)$ satisfying

(1)
\begin{equation}\label{2.2.10'}
 \rho\in C(R\times (0,T)), \ \ \rho(x,t)\ge 0, \ \ (x,t)\in R\times
(0,T);
\end{equation}

(2)
\begin{eqnarray}
 &&\sup_{t\in[0,T]}\int_{R}|\rho-\bar{\rho}|^2dx+\max_{(x,t)\in
R\times [0,T]}\rho\leq C; \label{2.2.11'}
\end{eqnarray}

(3)
\begin{equation}\label{2.2.12'}
\begin{aligned}
&\sup_{t\in[0,T]}\int_{R}(|\sqrt{\rho}u|^2+{(\rho^{\alpha-1/2})_x}^2+\frac{1}{\gamma-1}(\rho^\gamma-(\bar{\rho})^\gamma
-\gamma(\bar{\rho})^{\gamma-1}(\rho-\bar{\rho})))dx\\[3mm]
&+\int_{0}^{T}\int_{R}([(\rho^{\frac{\gamma+\alpha-1}{2}})_x]^2+{\Lambda(x,t)}^2)dxdt\leq
C,
\end{aligned}
\end{equation}
where $C$ is an absolute constant which only depends on the initial
data, and $\Lambda(x,t)\in L^2(R\times(0,T))$ is a function which
satisfies
$$
\int_{0}^{T}\int_{R}\Lambda\psi dxdt
=-\int_{0}^{T}\int_{R}{\rho}^{\alpha-1/2}\sqrt{\rho}u\psi_x
dxdt-\frac{2\alpha}{2\alpha-1}\int_0^T\int_R
{\rho}^{\alpha-1/2}_x\sqrt{\rho}u\psi dxdt.
$$
\end{theorem}

\begin{remark}
Under assumptions of  Theorem \ref{thm2.2.1'}, the case $\bar\rho=0$
has been dealt with in Proposition \ref{prop1.1}.
\end{remark}

The following is about the large time behavior of a weak solution.
\begin{theorem}\label{thm2.2.2}  Suppose that $(\rho(x,t),u(x,t))$ is a weak solution of the Cauchy problem $(\ref{2.2.2})-(\ref{2.2.4})$
satisfying $(\ref{2.2.10})-(\ref{2.2.16})$ or
$(\ref{2.2.10'})-(\ref{2.2.12'})$. Then
\begin{equation}\label{2.2.17}
\lim_{t\rightarrow +\infty}\sup_{x\in R}|\rho-\bar{\rho}|=0.
\end{equation}
\end{theorem}
\begin{remark}
In \cite{MV2008}, Mellet and Vasseur proved that if the initial data
is away from the vacuum  (has a positive lower bound) and
$0<\alpha<\frac 12$, the Cauchy problem \eqref{2.2.2}-\eqref{2.2.4}
has a unique global strong solution which is defined on $[0,T]$ for
any $T>0$. In comparison with  \cite{MV2008}, our results hold
uniform estimates  on $T$ and in the case that $\bar\rho=0$ the
vacuum at the infinity is permitted. Moreover,  by Theorem
\ref{thm2.2.2}, the large time behavior of the solutions of  the
strong solution can be obtained.
\end{remark}
Based on Theorem \ref{thm2.2.2}, it is easy to obtain
\begin{theorem}\label{thm2.2.3}
Suppose that the assumptions of Theorem \ref{thm2.2.1'} hold. Let $(\rho(x, t), u(x, t))$
be a weak solution of the Cauchy problem $(\ref{2.2.2})-(\ref{2.2.4})$ satisfying $(\ref{2.2.10'})-(\ref{2.2.12'})$.
Then for any $0 < \rho_1 < \bar\rho$, there exists a time $T_0$ such that
\begin{equation}
0 < \rho_1 \le\rho(x, t)\le C,\ \  (x, t)\in\ R \times[T_0,+\infty),
\end{equation}
where C is a constant same as in (\ref{2.2.11'}). Moreover, for $t \ge T_0$, the weak solution
becomes a unique strong solution to $(\ref{2.2.2})-(\ref{2.2.4})$, satisfying
$$\rho-\bar\rho\in L^\infty(T_0, t;H^1(R)),\ \  \rho_t \in L^\infty(T_0, t;L^2(R)),$$
$$u \in L^2(T_0, t;H^2(R)),\ \  u_t\in L^2(T_0, t;L^2(R))$$
and
\begin{equation}
\sup_{x\in R} |\rho-\bar\rho| + \|\rho-\bar\rho\|_{L^p(R)} + \|u\|_{L^2(R)}\rightarrow 0,
\end{equation}
as $t\rightarrow +\infty$, where $2 < p\le +\infty$.
\end{theorem}
\begin{remark}
Theorem \ref{thm2.2.3} shows that if $\bar\rho>0$, the vacuum will
vanish in finite time and the weak solution will become the strong
one after that. Similar to \cite{JX2008, LLX2008}, we can obtain
some results on the blow-up phenomena of the solutions when the
vacuum states vanish, which can be referred to \cite{JX2008,
LLX2008} for more details.
\end{remark}
\section{A Priori Estimates}

In this section, we will  construct approximate solutions and obtain
a priori estimates of the approximate solutions to the Cauchy
problem \eqref{2.2.2}-\eqref{2.2.4}. Two cases will be considered
respectively: $0<\alpha<\frac{1}{2}$ and $\alpha\ge\frac{1}{2}.$

{ \it Case I. $0<\alpha<\frac{1}{2}.$}

\vspace{3mm}

 For any given
$M>0$, we construct the smooth approximation solution of
$(\ref{2.2.2})-\eqref{2.2.4}$ on the cutoff domain
$\Omega^{M}=\{x\in R|-M<x<M\}$. Consider the initial condition
\begin{equation}\label{2.2.20}
(\rho, \rho u)(x, 0)=(\rho_{0\epsilon}, m_{0\epsilon}),
\end{equation}
and the boundary condition
\begin{equation}\label{2.2.21}
u(x,t)|_{x=\pm M}=0,
\end{equation}
where the initial data $\rho_{0\epsilon}, m_{0\epsilon}$ are smooth
functions satisfying
\begin{equation}\label{3.3}
\left\{\begin{array}{llll}
\rho_{0\epsilon}\rightarrow \rho_0\ {\rm in} \ L^1(\Omega^M)\cap L^\gamma(\Omega^M),\\
(\rho_{0\epsilon}^{\alpha-\frac{1}{2}})_x\rightarrow
(\rho_{0}^{\alpha-\frac{1}{2}})_x\ {\rm in} \ L^2(\Omega^M)\ \ {\rm
if}\  0<\alpha<\frac{1}{2},
 \\
(\log \rho_{0\epsilon})_x\rightarrow (\log \rho_{0})_x\ {\rm in} \ L^2(\Omega^M)\ \ {\rm if}\ \alpha=\frac{1}{2},\\
(m_{0\epsilon})^2(\rho_{0\epsilon})^{-1}\rightarrow (m_{0})^2(\rho_{0})^{-1},\ \ {\rm and} \\
(m_{0\epsilon})^{2+\delta}(\rho_{0\epsilon})^{-1-\delta}\rightarrow (m_{0})^{2+\delta}(\rho_{0})^{-1-\delta}\ {\rm in} \ L^1(\Omega^M),
\end{array}\right.
\end{equation}
as $\epsilon\rightarrow 0$. Here $\delta>0$, and there exists a
constant $C_0$ which does not depend on $\epsilon$ such that
\begin{equation}\label{3.4}
\rho_{0\epsilon}\geq C_0 \epsilon^{1/(2\alpha-2\theta)}.
\end{equation}

We note that the initial data can be regularized in an usual way(see
\cite{JWX2011} for instance).

The following estimate is a key one to prove the main theorem which
is based on the  energy and entropy estimates.
\begin{lemma}\label{lm2.1.1-}  Let
\begin{equation}\label{2.2.22}
\gamma>1, \ 0<\alpha<1/2.
\end{equation}
Assume that $(\rho_\epsilon, u_\epsilon)$ is the smooth solution of
$(\ref{2.2.2})$ with $\rho_{\epsilon}>0$. Then for any $T>0$, the
following estimate holds:
\begin{equation}\label{2.2.23}
\begin{aligned}
&\sup_{t\in [0,T]}{\int_{R}}\{\rho_{\epsilon}|u_{\epsilon}|^2+[(\frac{\rho_{\epsilon}^{\alpha-1/2}}{\alpha-{1/2}})_x]^2
+\rho_{\epsilon} \Psi(\rho_{\epsilon},\bar{\rho})\}(x,t)dx\\
&+\int_0^T\int_{R}\{\rho_{\epsilon}^{\alpha}(u_{\epsilon})_x^2+[(\rho_{\epsilon}^{\frac{\alpha+\gamma-1}{2}}
-{\bar{\rho}}^{\frac{\alpha+\gamma-1}{2}})_x]^2\}(x,t)dxdt\leq C,
\end{aligned}
\end{equation}
where $C$ is an universal constant independent of $\epsilon\ and \ T$.
\end{lemma}
\begin{proof}
It follows from $(\ref{2.2.2})_2$ that
\begin{equation}\label{2.2.24}
{\rho_\epsilon} {u_\epsilon}_t+{\rho_\epsilon} {u_\epsilon}{u_\epsilon}_x+({\rho_\epsilon}^\gamma)_x=({\rho_\epsilon}^\alpha {u_\epsilon}_x)_x.
\end{equation}
Multiply (\ref{2.2.24}) by ${u_\epsilon}$ to get
\begin{equation}\label{2.2.25}
(\frac{{\rho_\epsilon} |{u_\epsilon}|^2}{2})_t+(\frac{{\rho_\epsilon} {u_\epsilon}^3}{2})_x+{\rho_\epsilon}^\alpha({u_\epsilon}_x)^2+(P({\rho_\epsilon}))_x {u_\epsilon}-({\rho_\epsilon}^\alpha {u_\epsilon}{u_\epsilon}_x)_x=0
\end{equation}
In view of $(\ref{2.2.6})$, we have
\begin{equation}\label{2.2.26}
({\rho_\epsilon}\Psi({\rho_\epsilon},\bar{\rho}))_t+({\rho_\epsilon} {u_\epsilon}\Psi({\rho_\epsilon},\bar{\rho}))_x+{u_\epsilon}_x(P({\rho_\epsilon})-P(\bar{\rho}))=0.
\end{equation}
It follows from (\ref{2.2.25}) and (\ref{2.2.26}) that
\begin{equation}\label{2.2.27}
(\frac{{\rho_\epsilon} |{u_\epsilon}|^2}{2}+{\rho_\epsilon}\Psi({\rho_\epsilon},\bar{\rho}))_t+H_{1x}+{\rho_\epsilon}^\alpha({u_\epsilon}_x)^2=0,
\end{equation}
where $H_1=\frac{{\rho_\epsilon} {u_\epsilon}^3}{2}+{\rho_\epsilon} {u_\epsilon}\Psi({\rho_\epsilon},\bar{\rho})+{u_\epsilon}(P({\rho_\epsilon})-P(\bar{\rho}))-{\rho_\epsilon}^\alpha {u_\epsilon}{u_\epsilon}_x$.\\

Since
\begin{equation}\label{2.2.28}
({\rho_\epsilon}^\alpha {u_\epsilon}_x)_x=-{\rho_\epsilon} (\frac{{\rho_\epsilon}^{\alpha-1}}{\alpha-1})_{xt}-{\rho_\epsilon} {u_\epsilon}(\frac{{\rho_\epsilon}^{\alpha-1}}{\alpha-1})_{xx},
\end{equation}
 $(\ref{2.2.24})$ can be rewritten as
\begin{equation}\label{2.2.29}
{\rho_\epsilon} {u_\epsilon}_t+{\rho_\epsilon} {u_\epsilon}{u_\epsilon}_x+({\rho_\epsilon}^\gamma)_x=-{\rho_\epsilon} (\frac{{\rho_\epsilon}^{\alpha-1}}{\alpha-1})_{xt}-{\rho_\epsilon} {u_\epsilon}(\frac{{\rho_\epsilon}^{\alpha-1}}{\alpha-1})_{xx}.
\end{equation}
Multiplying (\ref{2.2.29}) by $(\frac{{\rho_\epsilon}^{\alpha-1}}{\alpha-1})_{x}$, we have
\begin{equation}\label{2.2.30}
\begin{aligned}
&(\frac{{\rho_\epsilon} (\frac{{\rho_\epsilon}^{\alpha-1}}{\alpha-1})_{x}^2}{2})_t+(\frac{{\rho_\epsilon} {u_\epsilon}(\frac{{\rho_\epsilon}^{\alpha-1}}{\alpha-1})_{x}^2}{2})_x
+({\rho_\epsilon} {u_\epsilon}(\frac{{\rho_\epsilon}^{\alpha-1}}{\alpha-1})_{x})_t+({\rho_\epsilon} {u_\epsilon}^2(\frac{{\rho_\epsilon}^{\alpha-1}}{\alpha-1})_{x})_x\\
&-{u_\epsilon}({\rho_\epsilon} (\frac{{\rho_\epsilon}^{\alpha-1}}{\alpha-1})_{xt}+{\rho_\epsilon} {u_\epsilon}(\frac{{\rho_\epsilon}^{\alpha-1}}{\alpha-1})_{xx})+(\frac{{\rho_\epsilon}^{\alpha-1}}{\alpha-1})_{x}(P({\rho_\epsilon}))_x=0.
\end{aligned}
\end{equation}
Multiplying $(\ref{2.2.29})$ by ${u_\epsilon}$ and adding up to
$(\ref{2.2.30})$, we obtain
\begin{equation}\label{2.2.31}
\begin{aligned}
&\{\frac{1}{2}{\rho_\epsilon}[{u_\epsilon}+(\frac{{\rho_\epsilon}^{\alpha-1}}{\alpha-1})_{x}]^2\}_t+\{\frac{1}{2}{\rho_\epsilon} {u_\epsilon}[{u_\epsilon}+(\frac{{\rho_\epsilon}^{\alpha-1}}{\alpha-1})_{x}]^2\}_x
+{u_\epsilon}(P({\rho_\epsilon}))_x\\
&+(\frac{{\rho_\epsilon}^{\alpha-1}}{\alpha-1})_{x}(P({\rho_\epsilon}))_x=0.
\end{aligned}
\end{equation}

From ($\ref{2.2.26}$) and $(\ref{2.2.31})$, we get
\begin{equation}\label{2.2.32}
\begin{aligned}
&\{\frac{1}{2}{\rho_\epsilon}[{u_\epsilon}+(\frac{{\rho_\epsilon}^{\alpha-1}}{\alpha-1})_{x}]^2+{\rho_\epsilon}\Psi({\rho_\epsilon},\bar{\rho})\}_t
+\{\frac{1}{2}{\rho_\epsilon} {u_\epsilon}[{u_\epsilon}+(\frac{{\rho_\epsilon}^{\alpha-1}}{\alpha-1})_{x}]^2+{\rho_\epsilon} {u_\epsilon} \Psi({\rho_\epsilon},\bar{\rho})\\
&+{u_\epsilon}(P({\rho_\epsilon})-P(\bar{\rho}))\}_x
+(\frac{{\rho_\epsilon}^{\alpha-1}}{\alpha-1})_{x}(P({\rho_\epsilon}))_x=0.
\end{aligned}
\end{equation}
Now we deal with the last term of the left hand side of
($\ref{2.2.32}$). Since
\begin{equation}\label{2.2.33}
(\frac{{\rho_\epsilon}^{\alpha-1}}{\alpha-1})_{x}={\rho_\epsilon}^{\alpha-2}{\rho_\epsilon}_x,
\end{equation}
we have
\begin{equation}\label{2.2.34}
(\frac{{\rho_\epsilon}^{\alpha-1}}{\alpha-1})_{x}(P({\rho_\epsilon}))_x={\rho_\epsilon}^{\alpha-2}{\rho_\epsilon}_x(P({\rho_\epsilon}))_x=
\frac{4\gamma}{(\gamma+\alpha-1)^2}[(\rho_{\epsilon}^{\frac{\alpha+\gamma-1}{2}}
-{\bar{\rho}}^{\frac{\alpha+\gamma-1}{2}})_x]^2.
\end{equation}
It follows from  ($\ref{2.2.32}$) and ($\ref{2.2.34}$) that
\begin{equation}\label{2.2.35}
\begin{aligned}
\{\frac{1}{2}{\rho_\epsilon}[{u_\epsilon}+(\frac{{\rho_\epsilon}^{\alpha-1}}{\alpha-1})_{x}]^2+{\rho_\epsilon}\Psi({\rho_\epsilon},\bar{\rho})\}_t
+H_{2x}+\frac{4\gamma}{(\gamma+\alpha-1)^2}[(\rho_{\epsilon}^{\frac{\alpha+\gamma-1}{2}}
-{\bar{\rho}}^{\frac{\alpha+\gamma-1}{2}})_x]^2=0,
\end{aligned}
\end{equation}
where
\begin{equation}\label{2.2.36}
H_2(x,t)=\frac{1}{2}{\rho_\epsilon} {u_\epsilon}[{u_\epsilon}+(\frac{{\rho_\epsilon}^{\alpha-1}}{\alpha-1})_{x}]^2+{\rho_\epsilon} {u_\epsilon} \Psi({\rho_\epsilon},\bar{\rho})+{u_\epsilon}(P({\rho_\epsilon})-P(\bar{\rho})).
\end{equation}
Multiplying ($\ref{2.2.36}$) by $\alpha$ and then adding up to ($\ref{2.2.27}$), we have
\begin{equation}\label{2.2.37}
\begin{aligned}
&\{\frac{\alpha}{2}{\rho_\epsilon}[{u_\epsilon}+(\frac{{\rho_\epsilon}^{\alpha-1}}{\alpha-1})_{x}]^2+\frac{{\rho_\epsilon} |{u_\epsilon}|^2}{2}+(\alpha+1){\rho_\epsilon}\Psi({\rho_\epsilon},\bar{\rho})\}_t
+[\alpha H_{2}+H_1]_x\\
&+{\rho_\epsilon}^\alpha {u_\epsilon}_x^2+\frac{4\gamma}{(\gamma+\alpha-1)^2}[(\rho_{\epsilon}^{\frac{\alpha+\gamma-1}{2}}
-{\bar{\rho}}^{\frac{\alpha+\gamma-1}{2}})_x]^2=0,
\end{aligned}
\end{equation}
Integrating ($\ref{2.2.37}$) over $[0, t] \times R$ with respect to $x,\ t$ gives
\begin{equation}\label{2.2.38}
\begin{aligned}
&\sup_{t\in [0,T]}\int_{R}\{{\rho_\epsilon} |{u_\epsilon}|^2+[(\frac{{\rho_\epsilon}^{\alpha-\frac{1}{2}}}{\alpha-\frac{1}{2}})_{x}]^2
+{\rho_\epsilon}\Psi({\rho_\epsilon},\bar{\rho})\}dx
+\int_0^T\int_R {\rho_\epsilon}^\alpha {u_\epsilon}_x^2dxdt\\
&+\int_0^T\int_R [(\rho_{\epsilon}^{\frac{\alpha+\gamma-1}{2}}
-{\bar{\rho}}^{\frac{\alpha+\gamma-1}{2}})_x]^2 dxdt\leq C.
\end{aligned}
\end{equation}
The proof of the lemma is finished.
\end{proof}
Based on Lemma \ref{lm2.1.1-}, we have
\begin{lemma}\label{lm2.1.2} Let
$0<\alpha<1/2$. Assume that $(\rho_\epsilon, u_\epsilon)$ is the
smooth solution of (\ref{2.2.2}) with $\rho_{\epsilon}>0$. Then
there exist two absolute  constants $C, \tilde C$ and a constant
$C(M)$ depending on $M$ such that
\begin{eqnarray}\label{2.2.39}
&&0<\tilde C\leq\rho_{\epsilon}\leq C,\ {\rm if}\ \bar{\rho}>0;\label{2.2.39}\\
&&0<C(M)\leq\rho_{\epsilon}\leq C,\ {\rm if}\
\bar{\rho}=0.\label{2.2.39+}
\end{eqnarray}
\end{lemma}
\begin{proof}\ \ In the case that  ${\bar{\rho}} > 0 $, from Lemma \ref{lm2.1.1}, we have
${\rho_\epsilon}\Psi({\rho_\epsilon},\bar{\rho})\in
L^{\infty}(0,T;L^1(R))$ and
$({{\rho_\epsilon}^{\alpha-\frac{1}{2}}})_x\in L^{\infty}(0,T;L^2$
$(R))$. Applying Lemma 5.3 in \cite{Lions1998}, we can get
\begin{equation}\label{2.2.43}
(\rho_\epsilon-\bar{\rho})1_{\{|\rho_\epsilon-\bar{\rho}|\le
\frac{\bar \rho}{2}\}}\in L^2(R)\ \ {\rm and}\
(\rho_\epsilon-\bar{\rho})1_{\{|\rho-_\epsilon\bar{\rho}|\ge
\frac{\bar\rho}{2}\}}\in L^\gamma(R).
\end{equation}
Since
\begin{equation}\label{2.2.43'}
{\rho_\epsilon}^\gamma-{\bar{\rho}}^\gamma-\gamma{\bar{\rho}}^{\gamma-1}({\rho_\epsilon}-\bar{\rho})\ge
C(\rho_\epsilon-\bar{\rho})^2 \ \ {\rm if}\ \gamma\ge2,
\end{equation}
thanks to \eqref{2.2.38}, we have
\begin{equation}\label{2.2.44}
\sup_{[0,T]}\int_R |{\rho_\epsilon}-\bar{\rho}|^2 dx\le C,\ \ {\rm if}\ \gamma\ge2,
\end{equation}
where $C$ is an universal constant independent of $\epsilon$ and
$T$.

 For $\bar\delta\in (0,\bar{\rho})$, we have
$|\rho_\epsilon|\leq\bar{\rho}+\bar\delta$ if
$|\rho_\epsilon-\bar{\rho}|\leq\bar\delta$, and hence \eqref{2.2.39}
holds true. If $|\rho_\epsilon-\bar{\rho}|\geq\bar\delta$, we can
prove that there exists a constant $C=C(\bar\delta)$ such that
\begin{equation}\label{2.2.45}
|\rho_{\epsilon}^s-{\bar{\rho}}^s|\leq C|\rho_\epsilon-\bar{\rho}|^s,\ \ \ s>0.
\end{equation}
In fact, from the facts
$$\frac{|\rho_{\epsilon}^s-{\bar{\rho}}^s|}{|\rho_\epsilon-\bar{\rho}|^s}\rightarrow 1,\ {\rm as}\ {\rho_\epsilon}\rightarrow\infty;\ \
\frac{|\rho_{\epsilon}^s-{\bar{\rho}}^s|}{|\rho_\epsilon-\bar{\rho}|^s}\rightarrow
1,\ {\rm as}\ {\rho_\epsilon}\rightarrow 0,$$
 there exist $\bar{\rho}_1,\bar{\rho}_2$ satisfying $0<\bar{\rho}_1<\bar{\rho}_2<\infty$ such that
$$
|\rho_{\epsilon}^s-{\bar{\rho}}^s|\leq 2|\rho_\epsilon-\bar{\rho}|^s,\ \ \rho_\epsilon\in[0,\bar{\rho}_1]\cup[\bar{\rho}_2,\infty).
$$

If $\rho_\epsilon\in[\bar{\rho}_1,\bar{\rho}_2],\
|\rho_\epsilon-\bar{\rho}|\geq \delta$, we have
$$
|\rho_{\epsilon}^s-{\bar{\rho}}^s|\leq
C|\rho_\epsilon-\bar{\rho}|^s,
$$
where $C$ depends on $\bar\delta, \bar{\rho}_1,\bar{\rho}_2$.  Thus
\eqref{2.2.45} holds true.\

It follows from (\ref{2.2.44})-(\ref{2.2.45}) that, for $\gamma\ge
2$,
\begin{equation}\label{2.2.48}
\begin{aligned}
|{\rho_\epsilon}-\bar{\rho}|^2&\leq\int_{R}({\rho_\epsilon}-\bar{\rho})^2dx+\int_{R}|2({\rho_\epsilon}-\bar{\rho})({\rho_\epsilon}-\bar{\rho})_x|dx\\
&=\int_{R}({\rho_\epsilon}-\bar{\rho})^2dx+\int_{R}|2({\rho_\epsilon}-\bar{\rho}){\rho_\epsilon}^{\frac{3}{2}-\alpha}
(\frac{{\rho_\epsilon}^{\alpha-\frac{1}{2}}}{\alpha-\frac{1}{2}})_x|dx\\
&\leq C+C(\int_{R}({\rho_\epsilon}-\bar{\rho})^2({\rho_\epsilon}^{\frac{3}{2}-\alpha}-{\bar{\rho}}^{\frac{3}{2}-\alpha})^2dx
+\int_{R}({\rho_\epsilon}-\bar{\rho})^2{\bar{\rho}}^{\frac{3}{2}-\alpha}dx)^{\frac{1}{2}}\\
&\leq C+C(\int_{R}({\rho_\epsilon}-\bar{\rho})^2({\rho_\epsilon}^{\frac{3}{2}-\alpha}-{\bar{\rho}}^{\frac{3}{2}-\alpha})^2 1|_{|
{\rho_\epsilon}-\bar{\rho}|\geq \bar\delta}dx\\
&\ \ \ \ \ \ \ \ \ \ \ \ \ +\int_{R}({\rho_\epsilon}-\bar{\rho})^2({\rho_\epsilon}^{\frac{3}{2}-\alpha}-{\bar{\rho}}^{\frac{3}{2}-\alpha})^2 1|_{|{\rho_\epsilon}-\bar{\rho}|\leq \bar \delta}dx)^{\frac{1}{2}}\\
&\leq C+C(\int_{R}|{\rho_\epsilon}-{\bar{\rho}}|^{5-2\alpha}1|_{|{\rho_\epsilon}-\bar{\rho}|\geq \delta}dx)^{\frac{1}{2}}\\
&\leq C+C \sup_{x\in R}|{\rho_\epsilon}-{\bar{\rho}}|^{\frac{3}{2}-\alpha}(\int_R ({\rho_\epsilon}-\bar{\rho})^2 dx)^{\frac{1}{2}}\\
&\leq C+C \sup_{x\in R}|{\rho_\epsilon}-{\bar{\rho}}|^{\frac{3}{2}-\alpha}.
\end{aligned}
\end{equation}
By Young's inequality and the condition $0<\alpha<\frac{1}{2}$, we
get
\begin{equation}\label{2.2.49}
|{\rho_\epsilon}-\bar{\rho}|^2\leq C,\ \ i.e.\ |{\rho_\epsilon}|\leq
C
\end{equation}
for $\gamma\ge 2$.

Now we consider the case $1<\gamma<2$. It follows from (\ref{2.2.43}) and (\ref{2.2.45}) that\\

\begin{eqnarray*}
&&|{\rho_\epsilon}-\bar{\rho}|^2
\leq\int_{\{|{\rho_\epsilon}-\bar{\rho}|>\frac{\bar{\rho}}{2}\}}|{\rho_\epsilon}-\bar{\rho}|^\gamma
dx
\sup_{x\in\{|{\rho_\epsilon}-\bar{\rho}|>\frac{\bar{\rho}}{2}\}}|{\rho_\epsilon}-\bar{\rho}|^{2-\gamma}\nonumber\\
&&+\int_{\{|{\rho_\epsilon}-\bar{\rho}|\le\frac{\bar{\rho}}{2}\}}|{\rho_\epsilon}-\bar{\rho}|^2
dx +\int_{R}|2({\rho_\epsilon}-\bar{\rho})({\rho_\epsilon}-\bar{\rho})_x|dx\nonumber\\
&&\le C \sup_{x\in
\{|{\rho_\epsilon}-\bar{\rho}|>\frac{\bar{\rho}}{2}\}}|{\rho_\epsilon}-\bar{\rho}|^{2-\gamma}
+C +\int_{R}
|2({\rho_\epsilon}-\bar{\rho}){\rho_\epsilon}^{\frac{3}{2}-\alpha}(\frac{{\rho_\epsilon}^{\alpha-\frac{1}{2}}}{\alpha-\frac{1}{2}})_x|dx\nonumber\\
&&\leq
C\sup_{\{|{\rho_\epsilon}-\bar{\rho}|>\frac{\bar{\rho}}{2}\}}|{\rho_\epsilon}-\bar{\rho}|^{2-\gamma}
+C+C[\int_{\{|\rho_\epsilon-\bar\rho|\ge
\frac{\bar\rho}{2}\}\cup\{|\rho_\epsilon-\bar\rho|\le
\frac{\bar\rho}{2}\}}
|{\rho_\epsilon}-\bar{\rho}|^2|{\rho_\epsilon}^{\frac{3}{2}-\alpha}-{\bar{\rho}}^{\frac{3}{2}-\alpha}|^2dx\nonumber\\
&& +\int_{\{|{\rho_\epsilon}-\bar{\rho}|>\frac{\bar{\rho}}{2}\}}({\rho_\epsilon}-\bar{\rho})^2({\bar{\rho}}^{\frac{3}{2}-\alpha})^2dx
+\int_{\{|{\rho_\epsilon}-\bar{\rho}|\le\frac{\bar{\rho}}{2}\}}({\rho_\epsilon}-\bar{\rho})^2({\bar{\rho}}^{\frac{3}{2}-\alpha})^2dx]^{\frac{1}{2}}\nonumber\\
\end{eqnarray*}
\begin{eqnarray}\label{2.2.48''.}
&&\leq
C\sup_{x\in\{|{\rho_\epsilon}-\bar{\rho}|>\frac{\bar{\rho}}{2}\}}|{\rho_\epsilon}-\bar{\rho}|^{2-\gamma}
+ C+C[\int_{\{|\rho_\epsilon-\bar\rho|\ge
\frac{\bar\rho}{2}\}}|{\rho_\epsilon}-\bar{\rho}|^{5-2\alpha}dx+\int_{\{|{\rho_\epsilon}-\bar{\rho}|>\frac{\bar{\rho}}{2}\}}
({\rho_\epsilon}-\bar{\rho})^2dx]^{\frac12}\nonumber\\
&&\leq
C+C\sup_{x\in\{|{\rho_\epsilon}-\bar{\rho}|>\frac{\bar{\rho}}{2}\}}|{\rho_\epsilon}-\bar{\rho}|^{2-\gamma}
+C[\sup_{x\in\{|{\rho_\epsilon}-\bar{\rho}|>\frac{\bar{\rho}}{2}\}}|{\rho_\epsilon}-\bar{\rho}|^{5-2\alpha-\gamma}+\sup_{x\in\{|{\rho_\epsilon}-\bar{\rho}|>\frac{\bar{\rho}}{2}\}}|{\rho_\epsilon}-\bar{\rho}|^{2-\gamma}]^{\frac12}\nonumber\\
&&\le C+C\sup_{x\in R}|{\rho_\epsilon}-\bar{\rho}|^{2-\gamma}
+C[\sup_{x\in R}|{\rho_\epsilon}-\bar{\rho}|^{5-2\alpha-\gamma}
+\sup_{x\in R}|{\rho_\epsilon}-\bar{\rho}|^{2-\gamma}]^{\frac12}.
\end{eqnarray}
By Young's inequality and the condition $0<\alpha<\frac{1}{2}$, we
get
\begin{equation}\label{2.2.49+}
|{\rho_\epsilon}-\bar{\rho}|^2\leq C,\ \ i.e.\ |{\rho_\epsilon}|\leq
C
\end{equation}
for $1<\gamma<2$. \eqref{2.2.49} and \eqref{2.2.49+} yield the
uniform upper bound in \eqref{2.2.39}.

 Now we prove the positive
lower bound estimate of $\rho_\epsilon$ in \eqref{2.2.39}.  Noticing
that $\lim_{{\rho_\epsilon}\rightarrow
0}{\rho_\epsilon}\Psi({\rho_\epsilon},\bar{\rho})=(\bar{\rho})^{\gamma}$,
we can obtain that ${\rho_\epsilon}\Psi({\rho_\epsilon},\bar{\rho})$
has a positive lower bound on $[0,\frac{1}{2}\bar{\rho}]$. Since
${\rho_\epsilon}\Psi({\rho_\epsilon},\bar{\rho})$ is bounded in
$L_T^{\infty}(L^1)$, there exists $C_1=C_1(T)>0$ such that for all
$t\in [0,T]$,
\begin{equation}\label{2.2.50}
{\rm meas}\{x\in R|\rho_\epsilon(x,t)\leq\frac{1}{2}\bar{\rho}\}\leq\frac{1}{\inf_{{\rho_\epsilon}\in[0,\frac{1}{2}\bar{\rho}]}{\rho_\epsilon}\Psi({\rho_\epsilon},\bar{\rho})}
\int_{\{x\in R|\rho_\epsilon\leq\frac{1}{2}\bar{\rho}\}}{\rho_\epsilon}\Psi({\rho_\epsilon},\bar{\rho})dx\leq C_1.
\end{equation}
Therefore for each $x_0\in R$, there exists $N=N(T)>0$ big enough
such that
$$\int_{|x-x_0|\leq N}\rho_{\epsilon}(x,t)dx
\geq \int_{\{|x-x_0|\leq N \}\cap\{x\in R|\rho_{\epsilon} >\frac{1}{2}\bar{\rho}\}}\rho_\epsilon dx$$
$$\geq\frac{1}{2}\bar{\rho} {\rm meas}\{\{|x-x_0|\leq N\}\cap\{x\in R|\rho_\epsilon(x,t)>\frac{1}{2}\bar{\rho}\}\}$$
\begin{equation}\label{2.2.51}
\begin{aligned}
&=\frac{1}{2}\bar{\rho} {\rm meas}\{\{|x-x_0|\leq N\}\cap\{x\in R|\rho_\epsilon(x,t)\leq\frac{1}{2}\bar{\rho}\}^c\}\\
&\geq \frac{1}{2}\bar{\rho}(2N-C_1)>0,\ \ \ \ t\in[0,T].
\end{aligned}
\end{equation}
Using the continuity of $\rho_\epsilon$, there exists
$x_1\in[x_0-N,x_0+N]$ such that
\begin{equation}\label{2.2.52}
\rho_\epsilon (x_1,t)=\int_{|x-x_0|\leq M}\rho_\epsilon (x,t)dx\geq \frac{1}{2}\bar{\rho}(2N-C_1).
\end{equation}
Then it follows from Lemma \ref{lm2.1.1} that
\begin{equation}\label{2.2.53}
\begin{aligned}
\rho_\epsilon^{\alpha-\frac{1}{2}}(x_0,t)
&=\rho_\epsilon^{\alpha-\frac{1}{2}}(x_1,t)+\int_{x_1}^{x_0}(\rho_\epsilon^{\alpha-\frac{1}{2}})_x(x,t)dx\\
&\leq [\frac{1}{2}\bar{\rho}(2N-C_1)]^{\alpha-\frac{1}{2}}+\|(\rho_\epsilon^{\alpha-\frac{1}{2}})_x(x,t)\|_{L^2(R)}|x_1-x_0|^{\frac{1}{2}}\\
&\leq [\frac{1}{2}\bar{\rho}(2N-C_1)]^{\alpha-\frac{1}{2}}+CN^{\frac{1}{2}}.
\end{aligned}
\end{equation}
Since $0<\alpha<\frac{1}{2}$, for any $x_0\in R$ and all
$t\in[0,T]$, we have
\begin{equation}\label{2.2.54}
\rho_\epsilon(x_0,t)
\geq \{[\frac{1}{2}\bar{\rho}(2N-C_1)]^{\alpha-\frac{1}{2}}+CN^{\frac{1}{2}}\}^{\frac{2}{2\alpha-1}}:=C(T).
\end{equation}
Up to now, we have proved  \eqref{2.2.39}.

 To prove \eqref{2.2.39+},  by
Lemma {\ref{lm2.1.1}}, we have
\begin{equation}\label{m1}
\begin{aligned}
{\rho_\epsilon}^\gamma&\leq\int_{R}{\rho_\epsilon}^\gamma dx+\int_{R}|2{\rho_\epsilon}^{\gamma-1}{\rho_\epsilon}_{x}|dx
\leq C+\int_{R}|2{\rho_\epsilon}^{\gamma-1}{\rho_\epsilon}^{\frac{3}{2}-\alpha}{\rho_\epsilon}^{\alpha-\frac{3}{2}}{\rho_\epsilon}_{x}|dx\\
&\leq C+\int_{R}|2{\rho_\epsilon}^{\gamma+\frac{1}{2}-\alpha}(\frac{{\rho_\epsilon}^{\alpha-\frac{1}{2}}}{\alpha-\frac{1}{2}})_x| dx
\leq C+(\int_R {\rho_\epsilon}^{2\gamma+1-2\alpha}dx)^{\frac{1}{2}}(\int_R ({\rho_\epsilon}^{\alpha-\frac{1}{2}})_x^2 dx)^{\frac{1}{2}}\\
&\leq C+C(\int_R {\rho_\epsilon}^\gamma dx)^{\frac{1}{2}}\sup_{x\in R}{\rho_\epsilon}^{\frac{\gamma+1-2\alpha}{2}}
\leq C+C\sup_{x\in R} {\rho_\epsilon}^{\frac{\gamma+1-2\alpha}{2}}.
\end{aligned}
\end{equation}
Applying  Young's inequality and the condition
$0<\alpha<\frac{1}{2}$, we obtain ${\rho_\epsilon}\leq C.$

To get the lower bound of $\rho_\epsilon$ in \eqref{2.2.39+}, we use
the Lagrangian coordinates as follows:
$$\xi=\int_{-M}^x
{\rho_\epsilon}(y, t)dy\ \quad \tau=t$$ where $x\in [-M, M], \
t>0$, and $\xi\in \Omega_L=(0, L)=(0, \int_{-M}^M
{\rho_\epsilon}(y,t)dy)=(0, \int_{-M}^M {\rho_\epsilon}_0(y)dy)$.
 In view of  the Lagrangian
coordinates transformation we get $({\rho_\epsilon}^\alpha)_{\xi}\in
L^2(\Omega_L)$ from Lemma \ref{lm2.1.1}. Let
$v=\frac{1}{{\rho_\epsilon}}$. Then we have
\begin{equation}\label{m2}
\begin{aligned}
v &\leq \int_{\Omega_L}v d\xi+\int_{\Omega_L} v^2|{\rho_\epsilon}_{\xi}|d\xi
\leq 2M+\frac{1}{2}\max_{\xi\in\Omega_L} v+\|({\rho_\epsilon}^\alpha)_\xi\|_{L^1(\Omega_L)}^{\frac{2}{1-2\alpha}}\\
&\leq2M+C+\frac{1}{2}\max_{\xi\in\Omega_L} v,
\end{aligned}
\end{equation}
which implies ${\rho_\epsilon}\geq C(M)>0$. \eqref{2.2.39+} is
proved and the proof of the lemma is finished.
\end{proof}

{\it Case II. $\alpha\ge \frac{1}{2}$.}

\vspace{3mm}

In the case that $\alpha\ge \frac{1}{2}$, we construct the
approximate solutions by solving
\begin{equation}\label{2.2.2''}
\left\{\begin{array}{llll}
 \rho_t+(\rho u)_x = 0 ,\\
 (\rho u)_t+(\rho u^2 +P(\rho))_x=(\mu_\epsilon(\rho) u_x)_x,\\
 (\rho_\epsilon,m_\epsilon)(x,t=0)=(\rho_\epsilon^0,m_\epsilon^0)
 \end{array}\right.
\end{equation}
with $\mu_\epsilon(\rho)=\rho^{\alpha}+\epsilon\rho^\theta,\
\epsilon>0,\ \theta\in (0,\frac12)$. The initial values
$(\rho_\epsilon^0, m_\epsilon^0)$ are regularized in the same way as
in \eqref{3.3} and \eqref{3.4}.

For any fixed $T > 0$ and for any fixed $\epsilon > 0$, there exists
a unique smooth approximate solution to (\ref{2.2.2''}) in the
region $(x, t) \in R \times(0, T)$.We refer to \cite{MV2008} for the
wellposedness of the global strong solution to the approximate
system (\ref{2.2.2''}).

Then we have
\begin{lemma}\label{lm2.1.1}  Let
\begin{equation}\label{2.2.22}
\gamma>1,\quad \alpha>1/2.
\end{equation}
Assume that $(\rho_\epsilon, u_\epsilon)$ is the smooth solution of
$(\ref{2.2.2})$ with $\rho_{\epsilon}>0$. Then for any $T>0$, the
following estimate holds:
\begin{equation}\label{2.2.23'}
\begin{aligned}
&\sup_{t\in [0,T]}{\int_{R}}\{\rho_{\epsilon}|u_{\epsilon}|^2+[(\frac{\rho_{\epsilon}^{\alpha-\frac12}}{\alpha-{\frac12}})_x]^2
+\epsilon^2[(\frac{\rho_{\epsilon}^{\theta-\frac12}}{\theta-{\frac12}})_x]^2+\rho_{\epsilon} \Psi(\rho_{\epsilon},\bar{\rho})\}(x,t)dx\\
&+\int_0^T\int_{R}\{(\rho_{\epsilon}^{\alpha}+\epsilon\rho_{\epsilon}^{\theta})[(u_{\epsilon})_x]^2+[(\rho_{\epsilon}^{\frac{\alpha+\gamma-1}{2}}
-{\bar{\rho}}^{\frac{\alpha+\gamma-1}{2}})_x]^2\\
&+[(\rho_{\epsilon}^{\frac{\theta+\gamma-1}{2}}
-{\bar{\rho}}^{\frac{\theta+\gamma-1}{2}})_x]^2\}(x,t)dxdt\leq C,
\end{aligned}
\end{equation}
where $C$ is an universal constant independent of $\epsilon\ and \ T$.
\end{lemma}

The proof is  similar to that of Lemma \ref{lm2.1.1-} (see also
Lemma 3.6 in \cite{JX2008}) and we omit the details here.

Based on this lemma, we have
\begin{lemma}\label{lm2.1.2'} Let $\alpha,\ \gamma$ satisfy (\ref{2.2.22}) and $\bar{\rho}>0$. Assume that $(\rho_\epsilon, u_\epsilon)$ is the smooth
solution of (\ref{2.2.2}) with $\rho_{\epsilon}>0$. Then there exist
an absolute constant $C$ and a constant $C(\epsilon,T)$ such that
\begin{eqnarray}\label{2.2.39'}
0<C(\epsilon,T)\leq\rho_{\epsilon}\leq C.
\end{eqnarray}
\end{lemma}
\begin{remark}
If $\bar\rho=0$, under the assumption \eqref{2.2.22},  the estimate
\eqref{2.2.39'} has been proved in \cite{JX2008}. If $\bar\rho>0$,
the estimate \eqref{2.2.39'} has also been proved in \cite{JX2008}
under the restriction \eqref{prop1.2-2}.
\end{remark}

\begin{proof}[Proof of Lemma \ref{lm2.1.2'}]\ \ We first prove the upper bound for $\rho_\epsilon(x,t)$.
The proof is divided into the following cases.

\vspace{2mm}

 If $\frac12<\alpha\le\frac{\gamma+1}{2}$, it follows from
(\ref{2.2.43}), (\ref{2.2.45}) and the entropy estimate
(\ref{2.2.23'}) that
\begin{equation}\label{2.2.48.}
\begin{aligned}
&|{\rho_\epsilon^{\alpha-\frac12}}-\bar{\rho}^{\alpha-\frac12}|^{2}
\leq\int_{-\infty}^x(|{\rho_\epsilon^{\alpha-\frac12}}-\bar{\rho}^{\alpha-\frac12}|^{2})_xdx\\
\le&\int_{-\infty}^x2({\rho_\epsilon^{\alpha-\frac12}}-\bar{\rho}^{\alpha-\frac12})({\rho_\epsilon^{\alpha-\frac12}}-\bar{\rho}^{\alpha-\frac12})_x|dx\\
\le& \int_{R}[({\rho_\epsilon^{\alpha-\frac12}}-\bar{\rho}^{\alpha-\frac12})_x]^2dx
+\int_{R}|\rho_\epsilon^{\alpha-\frac12}-\bar{\rho}^{\alpha-\frac12}|^{2}dx\\
\leq& C+\int_{R}({\rho_\epsilon^{\alpha-\frac12}}-\bar{\rho}^{\alpha-\frac12})^21_{\{|\rho_\epsilon-\bar{\rho}|<\frac{\bar{\rho}}{2}\}}dx
+\int_{R}({\rho_\epsilon^{\alpha-\frac12}}-\bar{\rho}^{\alpha-\frac12})^21_{\{|\rho_\epsilon-\bar{\rho}|\ge\frac{\bar{\rho}}{2}\}}dx\\
\leq& C+ I_1+I_2.
\end{aligned}
\end{equation}
Note that when $|\rho_\epsilon-\bar\rho| < \frac{\bar\rho}{2}$, that
is, $\frac{\bar\rho}{2}<\rho_\epsilon<\frac{3\bar\rho}{2}$, one has
$$|{\rho_\epsilon^{\alpha-\frac12}}-\bar{\rho}^{\alpha-\frac12}|^2\le |\rho_\epsilon-\bar\rho|^2\leq\rho_\epsilon\Psi(\rho_\epsilon,\bar\rho).$$
Hence,
\begin{equation}\label{2.2.48.1}
I_1\leq\int_{R}\rho_\epsilon\Psi(\rho_\epsilon,\bar\rho)dx\leq C.
\end{equation}
It follows from (\ref{2.2.45})that
\begin{equation}\label{2.2.48.2}
I_2\leq\int_{\{|\rho_\epsilon-\bar{\rho}|\ge\frac{\bar{\rho}}{2}\}}|\rho_\epsilon-\bar\rho|^{2(\alpha-\frac12)}dx
\leq\int_{\{|\rho_\epsilon-\bar{\rho}|\ge\frac{\bar{\rho}}{2}\}}C|\rho_\epsilon-\bar\rho|^\gamma dx
\leq\int_{R}\rho_\epsilon\Psi(\rho_\epsilon,\bar\rho)dx\leq C,
\end{equation}
if $\frac12<\alpha\le\frac{\gamma+1}{2}$. In view of
(\ref{2.2.48.})-(\ref{2.2.48.2}), we get
\begin{equation}
|{\rho_\epsilon^{\alpha-\frac12}}-\bar{\rho}^{\alpha-\frac12}|^{2}\le
C.
\end{equation}
Therefore $\rho_\epsilon$ has upper bound in the case that
$\frac12<\alpha\le\frac{\gamma+1}{2}$.

\vspace{2mm}

 If $\alpha>\frac12,\ 1<\gamma<2\alpha-1$, the proof is divided into the following subcases.

When    $\alpha>\frac12$, $\gamma\ge2$, it is easy to get that
\begin{equation}\label{2.2.43''}
{\rho_\epsilon}^\gamma-{\bar{\rho}}^\gamma-\gamma{\bar{\rho}}^{\gamma-1}({\rho_\epsilon}-\bar{\rho})\ge C(\rho_\epsilon-\bar{\rho})^2
\ \ {\rm if}\ \rho\ge 0.
\end{equation}
It follows from \eqref{2.2.23'} that
\begin{equation}\label{2.2.44'}
\sup_{[0,T]}\int_R |{\rho_\epsilon}-\bar{\rho}|^2 dx\le C,\ \ {\rm if}\ \gamma\ge2,
\end{equation}
where $C$ is an universal constant independent of $\epsilon$ and
$T$. Using (\ref{2.2.45}) and (\ref{2.2.44'}), for $\gamma\ge2$, we
have
\begin{equation}\label{2.2.48'}
\begin{aligned}
&|{\rho_\epsilon^{\alpha-\frac12}}-\bar{\rho}^{\alpha-\frac12}|^{2}
\leq\int_{-\infty}^x(|{\rho_\epsilon^{\alpha-\frac12}}-\bar{\rho}^{\alpha-\frac12}|^{2})_xdx\\
\le&\int_{-\infty}^x2({\rho_\epsilon^{\alpha-\frac12}}-\bar{\rho}^{\alpha-\frac12})({\rho_\epsilon^{\alpha-\frac12}}-\bar{\rho}^{\alpha-\frac12})_x|dx
\le C+\int_{R}|\rho_\epsilon^{\alpha-\frac12}-\bar{\rho}^{\alpha-\frac12}|^{2}dx\\
\le& C+\sup_{x\in R}|\rho_\epsilon^{\alpha-\frac12}-\bar{\rho}^{\alpha-\frac12}|^{2-\frac{2}{\alpha-\frac12}}\int_{R}|\rho_\epsilon-\bar{\rho}|^2dx\\
\le& C+C\sup_{x\in R}|\rho_\epsilon^{\alpha-\frac12}-\bar{\rho}^{\alpha-\frac12}|^{2-\frac{2}{\alpha-\frac12}}.
\end{aligned}
\end{equation}
By Young's inequality and the condition $\alpha>\frac{1}{2}$, we get
\begin{equation}\label{2.2.49-1}
|{\rho_\epsilon^{\alpha-\frac12}}-\bar{\rho}^{\alpha-\frac12}|^{2}\leq C,\ \ i.e.\ |{\rho_\epsilon}|\leq C.
\end{equation}

When  $\frac12<\alpha<\frac32$, $1<\gamma<2$, it follows from
(\ref{2.2.43}), (\ref{2.2.45}) and (\ref{2.2.48''.}) that
\begin{equation}\label{2.2.48'.}
\begin{aligned}
&|{\rho_\epsilon}-\bar{\rho}|^2\\
\leq&\int_{\{|{\rho_\epsilon}-\bar{\rho}|>\frac{\bar{\rho}}{2}\}}|{\rho_\epsilon}-\bar{\rho}|^\gamma dx
\sup_{\{|{\rho_\epsilon}-\bar{\rho}|>\frac{\bar{\rho}}{2}\}}|{\rho_\epsilon}-\bar{\rho}|^{2-\gamma}
+\int_{\{|{\rho_\epsilon}-\bar{\rho}|\le\frac{\bar{\rho}}{2}\}}|{\rho_\epsilon}-\bar{\rho}|^2 dx\\
&\ \ \ \ \ \ \ \ \ \ \ \ \ \ \ \ \ \ \ \ \ +\int_{R}|2({\rho_\epsilon}-\bar{\rho})({\rho_\epsilon}-\bar{\rho})_x|dx\\
\le& C+C\sup_{x\in R}|{\rho_\epsilon}-\bar{\rho}|^{2-\gamma}
+C(\sup_{x\in R}|{\rho_\epsilon}-\bar{\rho}|^{5-2\alpha-\gamma}
+\sup_{x\in R}|{\rho_\epsilon}-\bar{\rho}|^{3-2\alpha}\\
&\ \ \ \ \ \ \ \ \ \ \ \ \ \ \ \ \ \ \ \ \ +\sup_{x\in R}|{\rho_\epsilon}-\bar{\rho}|^{2-\gamma})^{\frac12}.
\end{aligned}
\end{equation}
By Young's inequality and the condition $\frac12<\alpha<\frac{3}{2}$ and $1<\gamma<2$, we can obtain
\begin{equation}\label{2.2.49-2}
|{\rho_\epsilon}-\bar{\rho}|^2\leq C,\ \ i.e.\ |{\rho_\epsilon}|\leq C.
\end{equation}

When  $\alpha\ge\frac32$, $1<\gamma<2$, it follows from (\ref{2.2.43}) and (\ref{2.2.45}) that
\begin{eqnarray}\label{2.2.48''}
&&|{\rho_\epsilon^{\alpha-\frac12}}-\bar{\rho}^{\alpha-\frac12}|^{2}
\leq\int_{-\infty}^x(|{\rho_\epsilon^{\alpha-\frac12}}-\bar{\rho}^{\alpha-\frac12}|^{2})_xdx\nonumber\\
&&\le\int_{-\infty}^x|2({\rho_\epsilon^{\alpha-\frac12}}-\bar{\rho}^{\alpha-\frac12})({\rho_\epsilon^{\alpha-\frac12}}-\bar{\rho}^{\alpha-\frac12})_x|dx
\le
C+\int_{-\infty}^x({\rho_\epsilon^{\alpha-\frac12}}-\bar{\rho}^{\alpha-\frac12})^2dx\nonumber\\
&&\leq C+C(\int_{R}({\rho_\epsilon^{\alpha-\frac12}}-\bar{\rho}^{\alpha-\frac12})^21_{\{|\rho_\epsilon-\bar{\rho}|\le\frac{\bar{\rho}}{2}\}}dx
+\int_{R}({\rho_\epsilon^{\alpha-\frac12}}-\bar{\rho}^{\alpha-\frac12})^21_{\{|\rho_\epsilon-\bar{\rho}|\ge\frac{\bar{\rho}}{2}\}}dx\nonumber\\
&&\leq C+\sup_{x\in R}|{\rho_\epsilon}^{\alpha-\frac{1}{2}}-{\bar{\rho}}^{\alpha-\frac{1}{2}}|^{2-\frac{2}{\alpha-\frac{1}{2}}}
\int_{R}({\rho_\epsilon}-\bar{\rho})^2 1_{\{|{\rho_\epsilon}-\bar{\rho}|\leq \frac{\bar{\rho}}{2}\}}dx\nonumber\\
&&+\sup_{x\in R}|{\rho_\epsilon}^{\alpha-\frac{1}{2}}-{\bar{\rho}}^{\alpha-\frac{1}{2}}|^{2-\frac{\gamma}{\alpha-\frac{1}{2}}}
\int_{R}|{\rho_\epsilon}-\bar{\rho}|^\gamma 1_{\{|{\rho_\epsilon}-\bar{\rho}|\geq \frac{\bar{\rho}}{2}\}}dx\nonumber\\
&&\le C+C\sup_{x\in R}|{\rho_\epsilon}^{\alpha-\frac{1}{2}}-{\bar{\rho}}^{\alpha-\frac{1}{2}}|^{2-\frac{2}{\alpha-\frac{1}{2}}}
+C\sup_{x\in R}|{\rho_\epsilon}^{\alpha-\frac{1}{2}}-{\bar{\rho}}^{\alpha-\frac{1}{2}}|^{2-\frac{\gamma}{\alpha-\frac{1}{2}}}.
\end{eqnarray}
By Young's inequality and the condition $\alpha\ge\frac{3}{2}$ and $1<\gamma<2$, we get
\begin{equation}\label{2.2.49-3}
|{\rho_\epsilon^{\alpha-\frac12}}-\bar{\rho}^{\alpha-\frac12}|^{2}\leq C,\ \ {\rm i.\ e.}\ |{\rho_\epsilon}|\leq C.
\end{equation}
Combining \eqref{2.2.49-1}, \eqref{2.2.49-2} with \eqref{2.2.49-3},
we get the uniform upper bound of $\rho_\epsilon$ for any
$\alpha>\frac12$ and $\gamma>1$.

 Next, we  prove the positive
lower bound estimate of $\rho_\epsilon$.

 Noticing that
$\lim_{{\rho_\epsilon}\rightarrow
0}{\rho_\epsilon}\Psi({\rho_\epsilon},\bar{\rho})=(\bar{\rho})^{\gamma}$,
we can obtain that ${\rho_\epsilon}\Psi({\rho_\epsilon},\bar{\rho})$
has a positive lower bound on $[0,\frac{1}{2}\bar{\rho}]$. Since
${\rho_\epsilon}\Psi({\rho_\epsilon},\bar{\rho})$ is bounded in
$L_T^{\infty}(L^1)$, there exists $C_1=C_1(T)>0$ such that for all
$t\in [0,T]$,
\begin{equation}\label{2.2.50'}
{\rm meas}\{x\in R|\rho_\epsilon(x,t)\leq\frac{1}{2}\bar{\rho}\}\leq\frac{1}{\inf_{{\rho_\epsilon}\in[0,\frac{1}{2}\bar{\rho}]}{\rho_\epsilon}\Psi({\rho_\epsilon},\bar{\rho})}
\int_{\{x\in R|\rho_\epsilon\leq\frac{1}{2}\bar{\rho}\}}{\rho_\epsilon}\Psi({\rho_\epsilon},\bar{\rho})dx\leq C_1.
\end{equation}
For each $x_0\in R$, there exists $N=N(T)>0$ large enough such that
$$\int_{|x-x_0|\leq N}\rho_{\epsilon}(x,t)dx
\geq \int_{\{|x-x_0|\leq N \}\cap\{x\in R|\rho_{\epsilon} >\frac{1}{2}\bar{\rho}\}}\rho_\epsilon dx$$
$$\geq\frac{1}{2}\bar{\rho} {\rm meas}\{\{|x-x_0|\leq N\}\cap\{x\in R|\rho_\epsilon(x,t)>\frac{1}{2}\bar{\rho}\}\}$$
\begin{equation}\label{2.2.51'}
\begin{aligned}
&=\frac{1}{2}\bar{\rho} {\rm meas}\{\{|x-x_0|\leq N\}\cap\{x\in R|\rho_\epsilon(x,t)\leq\frac{1}{2}\bar{\rho}\}^c\}\\
&\geq \frac{1}{2}\bar{\rho}(2N-C_1)>0,\ \ \ \ t\in[0,T].
\end{aligned}
\end{equation}
As the continuity of $\rho_\epsilon$, there exists $x_1\in[x_0-N,x_0+N]$ such that
\begin{equation}\label{2.2.52'}
\rho_\epsilon (x_1,t)=\int_{|x-x_0|\leq M}\rho_\epsilon (x,t)dx\geq \frac{1}{2}\bar{\rho}(2N-C_1).
\end{equation}
Then we can get from Lemma \ref{lm2.1.2} that
\begin{equation}\label{2.2.53'}
\begin{aligned}
\rho_\epsilon^{\theta-\frac{1}{2}}(x_0,t)
&=\rho_\epsilon^{\theta-\frac{1}{2}}(x_1,t)+\int_{x_1}^{x_0}(\rho_\epsilon^{\theta-\frac{1}{2}})_x(x,t)dx\\
&\leq [\frac{1}{2}\bar{\rho}(2N-C_1)]^{\theta-\frac{1}{2}}+\|(\rho_\epsilon^{\theta-\frac{1}{2}})_x(x,t)\|_{L^2(R)}|x_1-x_0|^{\frac{1}{2}}\\
&\leq [\frac{1}{2}\bar{\rho}(2N-C_1)]^{\alpha-\frac{1}{2}}+C_\epsilon N^{\frac{1}{2}}.
\end{aligned}
\end{equation}
For $\alpha>\frac{1}{2}$ and $\gamma>1$, for any $x_0\in R$ and all $t\in[0,T]$, we have
\begin{equation}\label{2.2.54'}
\rho_\epsilon(x_0,t)
\geq \{[\frac{1}{2}\bar{\rho}(2N-C_1)]^{\alpha-\frac{1}{2}}+C_\epsilon N^{\frac{1}{2}}\}^{\frac{2}{2\alpha-1}}:=C(\epsilon,T).
\end{equation}
The proof of the lemma is complete.
\end{proof}

Similarly, when $\alpha=\frac12$, we can establish the following a
priori estimates:
\begin{lemma}\label{lm2.1.3} Let $\alpha=1/2.$
Assume that $(\rho_\epsilon, u_\epsilon)$ is smooth solution of
(\ref{2.2.2}) with $\rho_{\epsilon}>0$. Then for any $T>0$, the
following estimate holds.
\begin{equation}\label{2.2.56}
\begin{aligned}
&\sup_{t\in [0,T]}\int_{R}\{{\rho_\epsilon} |{u_\epsilon}|^2+[({\rm
log}
{\rho_\epsilon})_{x}]^2+\epsilon^2[(\frac{{\rho_\epsilon}^{\theta-\frac{1}{2}}}{\theta-\frac{1}{2}})_x]^2
+{\rho_\epsilon}\Psi({\rho_\epsilon},\bar{\rho})\}dx
+\int_0^T\int_R {\rho_\epsilon}^{\frac{1}{2}} {u_\epsilon}_x^2dxdt\\
&+\int_0^T\int_R [(\rho_{\epsilon}^{\frac{\gamma-{\frac{1}{2}}}{2}}
-{\bar{\rho}}^{\frac{\gamma-{\frac{1}{2}}}{2}})_x]^2
+\epsilon[({\rho_\epsilon}^{\frac{\theta+\gamma-1}{2}}-{\bar{\rho}}^{\frac{\theta+\gamma-1}{2}})_x]^2dxdt\leq
C.
\end{aligned}
\end{equation}
Here $C$ is an universal constant independent of $\epsilon$ and $T$.
\end{lemma}
\begin{lemma}\label{lm2.1.4} Let
$\alpha=1/2$. Assume that $(\rho_\epsilon, u_\epsilon)$ is smooth
solution of (\ref{2.2.2}) with $\rho_{\epsilon}>0$. Then there exist
an absolute constant $C$ and a positive constant $C(\epsilon, T)$
such that
\begin{equation}\label{2.2.72}
0<C(\epsilon, T)\leq\rho_{\epsilon}\leq C.
\end{equation}
\end{lemma}
\begin{proof}\ \ Similar to the proof of lemma \ref{lm2.1.2'}, we first prove the positive uniform upper bound of $\rho_{\epsilon}$.

If $\bar{\rho}>0$, it follows from (\ref{2.2.56}) that
\begin{equation}\label{2.2.73}
\begin{aligned}
&|{\rho_\epsilon}-\bar{\rho}|^2\leq\int_{R}({\rho_\epsilon}-\bar{\rho})^2dx+\int_{R}|2({\rho_\epsilon}-\bar{\rho})({\rho_\epsilon}-\bar{\rho})_x|dx\\
&=\int_{R}({\rho_\epsilon}-\bar{\rho})^2dx+\int_{R}|2({\rho_\epsilon}-\bar{\rho}){\rho_\epsilon}({\rho_\epsilon}-\bar{\rho})_x{\rho_\epsilon}^{-1}|dx\\
&=\int_{R}({\rho_\epsilon}-\bar{\rho})^2dx+\int_{R}|2({\rho_\epsilon}-\bar{\rho}){\rho_\epsilon}(\log {\rho_\epsilon})_x|dx\\
&\leq
C+C(\int_{R}({\rho_\epsilon}-\bar{\rho})^2{\rho_\epsilon}^2dx)^{\frac{1}{2}}
(\int_{R}(\log {\rho_\epsilon})_x^2dx)^{\frac{1}{2}}\\
&\leq
C+C(\int_{R}({\rho_\epsilon}-\bar{\rho})^2({\rho_\epsilon}-{\bar{\rho}})^2dx
+\int_{R}({\rho_\epsilon}-\bar{\rho})^2{\bar{\rho}}^2dx)^{\frac{1}{2}}\\
&\leq
C+C(\int_{R}({\rho_\epsilon}-\bar{\rho})^2({\rho_\epsilon}-{\bar{\rho}})^2dx
+C)^{\frac{1}{2}}\\
&\leq C+\frac{1}{2} \sup_{x\in R}|{\rho_\epsilon}-{\bar{\rho}}|^2.
\end{aligned}
\end{equation}
By using Young's Inequality, we obtain that
\begin{equation}\label{2.2.74}
|{\rho_\epsilon}-\bar{\rho}|^2\leq C,\ \ i.e.\
{\rho_\epsilon}\leq\bar{\rho}+C.
\end{equation}

If ${\bar{\rho}} = 0$, it follows from (\ref{2.2.56}) that
\begin{equation}\label{m1}
\begin{aligned}
{\rho_\epsilon}^\gamma&\leq\int_{R}{\rho_\epsilon}^\gamma dx+\int_{R}|2{\rho_\epsilon}^{\gamma-1}{\rho_\epsilon}_{x}|dx\leq C+\int_{R}|2{\rho_\epsilon}^{\gamma-1}{\rho_\epsilon}{\rho_\epsilon}^{-1}{\rho_\epsilon}_{x}|dx\\
&\leq C+\int_{R}|2{\rho_\epsilon}^{\gamma}(\log {{\rho_\epsilon}})_x| dx\leq C+(\int_R {\rho_\epsilon}^{2\gamma}dx)^{\frac{1}{2}}(\int_R (\log{{\rho_\epsilon}})_x^2)^{\frac{1}{2}}\\
&\leq C+C(\int_R {\rho_\epsilon}^\gamma dx)^{\frac{1}{2}} \sup_{x\in
R} {\rho_\epsilon}^{\frac{\gamma}{2}}\leq C+C(\sup_{x\in R}
{\rho_\epsilon}^\gamma)^{\frac{1}{2}}.
\end{aligned}
\end{equation}
By using Young's Inequality, we obtain that ${\rho_\epsilon}\leq C.$

We now prove the positive lower bound estimate of $\rho_\epsilon$.
Using ($\ref{2.2.50}$)-($\ref{2.2.52}$) and (\ref{2.2.56}), we
obtain
\begin{equation}\label{2.2.75}
\begin{aligned}
\rho_\epsilon^{\theta-\frac{1}{2}}(x_0,t)
&=\rho_\epsilon^{\theta-\frac{1}{2}}(x_1,t)+\int_{x_1}^{x_0}(\rho_\epsilon^{\theta-\frac{1}{2}})_x(x,t)dx\\
&\leq [\frac{1}{2}\bar{\rho}(2N-C_1)]^{-1}+max\rho_\epsilon^{-1} \|(\rho_\epsilon^{\theta-\frac{1}{2}})_x(x,t)\|_{L^2(R)}|x_1-x_0|^{\frac{1}{2}}\\
&\leq [\frac{1}{2}\bar{\rho}(2N-C_1)]^{-1}+C_\epsilon N.
\end{aligned}
\end{equation}
Since $0<\theta<\frac{1}{2}$, by the construction of the approximate
solutions in \eqref{2.2.2''}, we have
\begin{equation}\label{2.2.76}
\rho_\epsilon(x_0,t) \geq
C\{[\frac{1}{2}\bar{\rho}(2N-C_1)]^{-1}+C_\epsilon
N\}^{-1}:=C(\epsilon, T),
\end{equation}
for any $x_0\in R$ and $t\in[0,T]$.

The proof of the lemma is finished.
\end{proof}

\section{Proof of the Main Results }

In this section, we give the proof of the main results. The proof is
completely similar to those in \cite{GJX2008, JX2008, LLX2008} and
we give a sketch of proof here.

\begin{proof}[{Proof of Theorem \ref{thm2.2.1}.}] Based on a priori
estimates of Lemma \ref{lm2.1.1-}-Lemma \ref{lm2.1.2} and Lemma
\ref{lm2.1.3}-Lemma \ref{lm2.1.4}, applying similar approaches in
\cite{GJX2008, LLX2008, MV2007} and the references therein, we can
obtain that for any $T>0$ there exists a unique global smooth
solution of ($\ref{2.2.2}$)-($\ref{2.2.4}$) satisfying
$$\rho_{\epsilon}, \rho_{\epsilon x}, \rho_{\epsilon t}, u_{\epsilon}, u_{\epsilon x}, u_{\epsilon t},
u_{\epsilon xx}\in C^{\beta, \frac{\beta}{2}}([-M,M]\times[0, T]),\
\ 0<\beta<1,$$ and $\rho_\epsilon\geq C(\epsilon)>0$ in
$[-M,M]\times[0, T]$ when $0<\alpha\leq\frac{1}{2}$. And, the
estimates in Lemma \ref{lm2.1.1-}-Lemma \ref{lm2.1.2} and Lemma
\ref{lm2.1.3}-Lemma \ref{lm2.1.4} hold for
$\{\rho_\epsilon,u_\epsilon\}$.

We only give a  proof of the case $0<\alpha<\frac{1}{2}$. The case
$\alpha=\frac{1}{2}$ can be proved in a similar way. For any fixed
$M>0$, similar to \cite{GJX2008, JX2008, LLX2008}, we can obtain
that (up to a subsequence)
\begin{equation}\label{2.2.77}
\rho_\epsilon\rightarrow\rho \ \ \ \ in\ C([0,T]\times[-M,M]),
\end{equation}
\begin{equation}\label{2.2.78}
(\rho_\epsilon^{\alpha-\frac12})_x\rightharpoonup(\rho^{\alpha-\frac12})_x
\ \ \ weakly \ in \ L^2((0,T)\times[-M,M]),
\end{equation}
\begin{equation}\label{2.2.79}
\rho_\epsilon^{\alpha}u_{\epsilon x}\rightharpoonup \Lambda\ \ \ \
weakly\ in \ L^2((0,T)\times[-M,M]),
\end{equation}
as $\epsilon\rightarrow 0$, for some function $\rho\in
C([0,T]\times[-M,M])$ and $\Lambda\in L^2((0,T)\times[-M,M])$ which
satisfy
\begin{equation}\label{2.2.80}
\int_0^T\int_{-M}^M \Lambda\varphi
dxdt=-\int_0^T\int_{-M}^M\rho^{\alpha-\frac12}\sqrt{\rho}u\varphi_x
dxdt -\frac{2\alpha}{2\alpha-1}\int_0^T\int_{-M}^M
(\rho^{\alpha-\frac12})_x\sqrt{\rho}u\varphi dxdt.
\end{equation}
To get the convergence of the term $\sqrt{\rho_\epsilon}
u_\epsilon$, we apply similar approaches in \cite{GJX2008, JX2008,
LLX2008, MV2007}. More precisely, we have ${\rho_\epsilon}
u_\epsilon$ converges strongly in $L^1([0,T]\times[-M,M])$ and
$L^2([0,T];$ $L^{1+\zeta}(-M,M))$ and almost everywhere to some
function $m(x,t)$, where $\zeta>0$ is some small positive number.
Also, we can prove that $\sqrt{\rho_\epsilon} u_\epsilon$ converges
strongly in $L^2([0,T]\times[-M,M])$ to $\frac{m}{\sqrt{\rho}}$
which is defined to be zero when $m=0$ and there exists a function
$u(x,t)$ such that $m(x,t)=\rho(x,t) u(x,t)$. Moreover, we have
\begin{equation}\label{2.2.81}
\rho\in C(R\times (0,T)),
\end{equation}
\begin{equation}\label{2.2.82}
\sup_{t\in[0,T]}\int_{-M}^M|\rho-\bar{\rho}|^2dx+\max_{(x,t)\in
R\times [0,T]}\rho\leq C,
\end{equation}
\begin{equation}\label{2.2.83}
\begin{aligned}
\sup_{t\in[0,T]}\int_{-M}^M(|\sqrt{\rho}u|^2&+{(\rho^{\alpha-\frac12})_x}^2
+\frac{1}{\gamma-1}(\rho^\gamma-(\bar{\rho})^\gamma
-\gamma(\bar{\rho})^{\gamma-1}(\rho-\bar{\rho})))dx\\
+&\int_{0}^{T}\int_{-M}^M([(\rho^{\frac{\gamma+\alpha-1}{2}})_x]^2+{\Lambda(x,t)}^2)dxdt\leq
C,
\end{aligned}
\end{equation}
where C is an absolute constant depending on the initial data.

 Using
a diagonal procedure, we obtain that the above converges (up to a
subsequence) remain true for any $M>0$ and the existence of weak
solutions of (\ref{2.2.2})-(\ref{2.2.4}) can be directly proved.
Moreover, (\ref{2.2.10})-(\ref{2.2.16}) hold true due to
(\ref{2.2.81})-(\ref{2.2.83}). The proof of the theorem is complete.
\end{proof}

Using  Lemma \ref{lm2.1.1}-Lemma \ref{lm2.1.2'}, we can prove
Theorem \ref{thm2.2.1'}  in a similar way ( see also \cite{GJX2008,
LLX2008, MV2007}) and we omit the details here.

\vspace{3mm}

 Now we give the proof of Theorem  \ref{thm2.2.2}, which
is about the asymptotic behavior of the weak solutions. We assume
that the solutions are smooth enough. The rigorous proof can be
obtained by using the usual regularization procedure.

\begin{proof}[Proof of Theorem \ref{thm2.2.2}]

We only prove the case of $0<\alpha<\frac{1}{2}$ in Theorem
\ref{thm2.2.2}, since other cases can be proved in a similar way.
The proof can be done by considering the cases  $\bar{\rho}>0$ and
$\bar{\rho}=0$ respectively.

If $\bar{\rho}>0$, since $0\leq\rho\leq C,\ \bar{\rho}>0$, for some constant $C_1>0$, we have
\begin{equation}\label{2.2.84}
C_1^{-1}(\rho-\bar{\rho})^2\leq \rho\Psi(\rho,\bar{\rho})\leq C_1(\rho-\bar{\rho})^2.
\end{equation}
From  Lemma \ref{lm2.1.2}  we have $|\rho^s-{\bar{\rho}}^s|^2\leq
C|\rho-\bar{\rho}|^2$ for any $s>0$. Hence,
\begin{equation}\label{2.2.85}
\int_{R}|\rho^s-{\bar{\rho}}^s|^2dx\leq C\int_{R}|\rho-\bar{\rho}|^2dx\leq C.
\end{equation}
Similarly,
\begin{equation}\label{2.2.86}
\int_{R}|\rho^s-{\bar{\rho}}^s|^{2\lambda}dx\leq C\int_{R}|\rho-\bar{\rho}|^{2\lambda}dx\leq C,
\end{equation}
for any $s>0, \lambda>1$. Moreover, one has
\begin{equation}\label{2.2.87}
\begin{aligned}
&\int_{R}|[(\rho^s-{\bar{\rho}}^s)^{2\lambda}]_x|dx=2\lambda s\int_{R}|(\rho^s-{\bar{\rho}}^s)^{2\lambda-1}[\rho^{s-1}\rho_x]|dx\\
&\leq \frac{2\lambda s}{|\alpha-\frac{1}{2}|}
(\int_{R}(\rho^s-{\bar{\rho}}^s)^{2(2\lambda-1)}\rho^{2s+1-2\alpha}dx)^{\frac{1}{2}}(\int_{R}[(\rho^{\alpha-\frac{1}{2}})_x]^2dx)^{\frac{1}{2}}\\
&\leq C.
\end{aligned}
\end{equation}

By Lemma \ref{2.2.22}, one has
\begin{equation}\label{2.2.89}
\int_0^T\int_{R}[(\rho^{\frac{\alpha+\gamma-1}{2}}-{\bar{\rho}}^{\frac{\alpha+\gamma-1}{2}})_x]^2(x,t)dxdt\leq C
\end{equation}
Denote $b=\frac{\alpha+\gamma-1}{2}$. Then
\begin{equation}\label{2.2.90}
\int_0^T\int_{R}[(\rho^b-{\bar{\rho}}^b)_x]^2(x,t)dxdt\leq C
\end{equation}
Choosing $s>b+1$, one has
\begin{equation}\label{2.2.91}
\begin{aligned}
&(\rho^s-{\bar{\rho}}^s)^{2}=\int_{-\infty}^x[(\rho^s-{\bar{\rho}}^s)^2]_x dx
=2\int_{-\infty}^x (\rho^s-{\bar{\rho}}^s)(\rho^s-{\bar{\rho}}^s)_x dx\\
&=2s\int_{-\infty}^x (\rho^s-{\bar{\rho}}^s)(\rho^{s-1}{\rho}_x) dx
=\frac{2s}{b}\int_{-\infty}^x (\rho^s-{\bar{\rho}}^s)[(\rho^b-{\bar{\rho}}^b)_x \rho^{s-b}]dx\\
&\leq \|\rho^s-{\bar{\rho}}^s\|_{L^2(R)}\|(\rho^b-{\bar{\rho}}^b)_x\|_{L^2(R)}.
\end{aligned}
\end{equation}
Consequently,
\begin{equation}\label{2.2.92}
\begin{aligned}
\int_{0}^t\sup_{x\in R}(\rho^s-{\bar{\rho}}^s)^4 dt\leq C\sup_{x\in R}\|\rho^s-{\bar{\rho}}^s\|_{L^2(R)}^2
\int_{0}^t\|(\rho^b-{\bar{\rho}}^b)_x\|_{L^2(R)}^2 dt
\leq C.
\end{aligned}
\end{equation}
Moreover, applying($\ref{2.2.86}$), one has
\begin{equation}\label{2.2.93}
\begin{aligned}
&\int_{0}^t\int_{R}(\rho^s-{\bar{\rho}}^s)^4 (\rho^s-{\bar{\rho}}^s)^{2l}dxdt
\leq \int_{0}^t[\sup_{x\in R}(\rho^s-{\bar{\rho}}^s)^4\int_{R}(\rho^s-{\bar{\rho}}^s)^{2l}dx]dt\\
&\leq \sup_{t} \int_{R}(\rho^s-{\bar{\rho}}^s)^{2l}dx\int_{0}^t\sup_{x\in R}(\rho^s-{\bar{\rho}}^s)^4dt
\leq C,
\end{aligned}
\end{equation}
where $l\geq 1$ is any real number. Hence
\begin{equation}\label{2.2.94}
\int_{0}^t\int_{R}(\rho^s-{\bar{\rho}}^s)^{4+2l}dxdt\leq C.
\end{equation}
Denote $f(t)=\int_{R}(\rho^s-{\bar{\rho}}^s)^{4+2l}dx$. Then, from
($\ref{2.2.86}$) and ($\ref{2.2.94}$), one has $f(t)\in L^1(0,
\infty)\cap L^\infty(0, \infty)$. Furthermore, direct calculations
show that
\begin{equation}\label{2.2.95}
\begin{aligned}
&\frac{d}{dt}f(t)=(4+2l)s\int_{R}(\rho^s-{\bar{\rho}}^s)^{3+2l}\rho^{s-1}\rho_t dx\\
&=-(4+2l)s\int_{R}(\rho^s-{\bar{\rho}}^s)^{3+2l}\rho^{s-1}(\rho u)_x dx\\
&=(4+2l)(3+2l)s\int_{R}(\rho^s-{\bar{\rho}}^s)^{2+2l}(\rho^{s})_x\rho^{s-1}\rho u dx\\
&+(4+2l)s\int_{R}(\rho^s-{\bar{\rho}}^s)^{3+2l}(s-1)\rho^{s-2}\rho_x\rho u dx\\
&=(4+2l)(3+2l)s\int_{R}(\rho^s-{\bar{\rho}}^s)^{2+2l}(\rho^{s})_x\rho^{s-1}\rho u dx\\
&+(4+2l)s(s-1)\int_{R}(\rho^s-{\bar{\rho}}^s)^{3+2l}\rho^{s-2}\rho_x\rho u dx\\
&=J_1+J_2
\end{aligned}
\end{equation}
Now, we claim that $J_i(t)\in L^2(0,+\infty),\ (i=1,2)$. In fact,
\begin{equation}\label{2.2.96}
\begin{aligned}
J_1(t)&=(4+2l)(3+2l)s\int_{R}(\rho^s-{\bar{\rho}}^s)^{2+2l}(\rho^{s})_x\rho^{s-1}\rho u dx\\
&=\frac{(4+2l)(3+2l)s^2}{b}\int_{R}(\rho^s-{\bar{\rho}}^s)^{2+2l}(\rho^{b})_x\rho^{2s-b} u dx\\
&=\frac{(4+2l)(3+2l)s^2}{b}\int_{R}(\rho^s-{\bar{\rho}}^s)^{2+2l}(\rho^{b})_x\rho^{2s-b-\frac{1}{2}}\sqrt{\rho} u dx\\
&\leq C\|\sqrt{\rho} u\|_{L^2(R)}\|(\rho^{b}-{\bar{\rho}}^b)_x\|_{L^2(R)},
\end{aligned}
\end{equation}
Hence, by Lemma \ref{lm2.1.1},
\begin{equation}\label{2.2.97}
\int_0^t|J_1(t)|^2 dt\leq C\sup_{t\in[0,T]}\|\sqrt{\rho} u\|_{L^2(R)}^2\int_0^t\|(\rho^{b}-{\bar{\rho}}^b)_x\|_{L^2(R)}^2dt
\leq C.
\end{equation}
And
\begin{equation}\label{2.2.98}
\begin{aligned}
J_2(t)&=(4+2l)s(s-1)\int_{R}(\rho^s-{\bar{\rho}}^s)^{3+2l}\rho^{s-1}\rho_x u dx\\
&=(4+2l)(s-1)\int_{R}(\rho^s-{\bar{\rho}}^s)^{3+2l}({\rho^{s}})_x u dx\\
&=\frac{(4+2l)(s-1)s}{b}\int_{R}(\rho^s-{\bar{\rho}}^s)^{3+2l}({\rho^{b}})_x \rho^{s-b}u dx\\
&=\frac{(4+2l)(s-1)s}{b}\int_{R}(\rho^s-{\bar{\rho}}^s)^{3+2l}({\rho^{b}})_x \rho^{s-b-\frac{1}{2}}\sqrt{\rho}u dx\\
&\leq C\|\sqrt{\rho} u\|_{L^2(R)}\|(\rho^{b}-{\bar{\rho}}^b)_x\|_{L^2(R)},
\end{aligned}
\end{equation}
Using Lemma \ref{lm2.1.1} again, one has
\begin{equation}\label{2.2.99}
\int_0^t|J_2(t)|^2 dt\leq C\sup_{t\in[0,T]}\|\sqrt{\rho} u\|_{L^2(R)}^2\int_0^t\|(\rho^{b}-{\bar{\rho}}^b)_x\|_{L^2(R)}^2dt
\leq C.
\end{equation}
Consequently,
\begin{equation}\label{2.2.100}
\frac{d}{dt} f(t)\in L^2(0,+\infty).
\end{equation}
Combining the obtained fact that $f(t)\in L^1(0, \infty)\cap L^\infty(0, \infty)$, one has
\begin{equation}\label{2.2.101}
f(t)\rightarrow 0,\ \ t\rightarrow +\infty.
\end{equation}
Letting $m\geq 1$ be any real number to be determined later, we have
\begin{equation}\label{2.2.102}
\begin{aligned}
|\rho^s-{\bar{\rho}}^s|^{m}&=|\int_{-\infty}^x[(\rho^s-{\bar{\rho}}^s)^m]_x dx|
=|m\int_{-\infty}^x (\rho^s-{\bar{\rho}}^s)^{m-1}(\rho^s-{\bar{\rho}}^s)_x dx|\\
&=|m\int_{-\infty}^x (\rho^s-{\bar{\rho}}^s)^{m-1}[\frac{s}{\alpha-\frac{1}{2}}(\rho^{s-\alpha+\frac{1}{2}}({\rho^{\alpha-\frac{1}{2}}})_x)] dx|\\
&\leq C\|(\rho^s-{\bar{\rho}}^s)^{m-1}\|_{L^2(R)}\|(\rho^{\alpha-\frac{1}{2}})_x\|_{L^2(R)}\\
&\leq C\|(\rho^s-{\bar{\rho}}^s)^{m-1}\|_{L^2(R)}.
\end{aligned}
\end{equation}
Choosing $2(m-1)=4+2l$, one has
\begin{equation}\label{2.2.103}
\sup_{x\in R}|\rho^s-{\bar{\rho}}^s|^{m}\leq Cf^{\frac{1}{2}}(t)\rightarrow 0,\ \ t\rightarrow 0.
\end{equation}
Therefore, $\lim_{t\rightarrow+\infty}\sup_{x\in
R}|\rho^s-{\bar{\rho}}^s|=0$. Using the fact that
$$|\rho-{\bar{\rho}}|^s=|\rho-{\bar{\rho}}|^s1|_{0\leq\rho\leq\frac{\bar{\rho}}{2}}+|\rho-{\bar{\rho}}|^s1|_{\rho>\frac{\bar{\rho}}{2}}
\leq C|\rho^s-{\bar{\rho}}^s|1|_{0\leq\rho\leq\frac{\bar{\rho}}{2}}+C|\rho^s-{\bar{\rho}}^s|^s1|_{\rho>\frac{\bar{\rho}}{2}}.$$
Hence,
$$\sup_{x\in R}|\rho-{\bar{\rho}}|^s\leq C\sup_{x\in R}|\rho^s-{\bar{\rho}}^s|+C\sup_{x\in R}|\rho^s-{\bar{\rho}}^s|^s
\rightarrow 0,\ \ t\rightarrow +\infty,$$
which implies that
$$\lim_{t\rightarrow+\infty}\sup_{x\in R}|\rho-{\bar{\rho}}|=0.$$

If  $\bar{\rho}=0$, from  Lemma \ref{lm2.1.2}, we have that $\rho^{2s}\leq C\rho^\gamma$
for any $s>\frac{\gamma}{2}$. Hence,
\begin{equation}\label{2.2.85'}
\int_{R}\rho^{2s}dx\leq C\int_{R}\rho^\gamma dx\leq C.
\end{equation}
Similarly,
\begin{equation}\label{2.2.86'}
\int_{R}(\rho^s)^{2\lambda}dx\leq\int_{R}\rho^{2s\lambda}dx\leq C\int_{R}\rho^{\gamma\lambda}dx\leq C,
\end{equation}
for any $s>\frac{\gamma}{2}, \lambda>1$. Moreover, one has
\begin{equation}\label{2.2.87'}
\begin{aligned}
&\int_{R}|[(\rho^{2s\lambda}]_x|dx=2\lambda s\int_{R}|\rho^{2s\lambda-1}\rho_x|dx=2\lambda s\int_{R}|\rho^{2s\lambda-1}\rho^{\frac{3}{2}-\alpha}\rho^{\alpha-\frac{3}{2}}\rho_x|dx\\
\leq& \frac{2s\lambda}{|\alpha-\frac{1}{2}|}
(\int_{R}\rho^{4s\lambda-2\alpha+1}dx)^{\frac{1}{2}}(\int_{R}[(\rho^{\alpha-\frac{1}{2}})_x]^2dx)^{\frac{1}{2}}
\leq C.
\end{aligned}
\end{equation}

Denote $b=\frac{\alpha+\gamma-1}{2}$. Then we have
\begin{equation}\label{2.2.90'}
\int_0^T\int_{R}[(\rho^b)_x]^2(x,t)dxdt\leq C,
\end{equation}
by Lemma \ref{lm2.1.1} and $\bar{\rho}=0$.
Choosing $s>b+\frac{\gamma}{2}$, one has
\begin{equation}\label{2.2.91'}
\begin{aligned}
&\rho^{2s}=\int_{-\infty}^x(\rho^{2s})_x dx
=2s\int_{-\infty}^x \rho^{2s-1}\rho_x dx
=\frac{2s}{b}\int_{-\infty}^x \rho^{\frac{2s+1-(\alpha+\gamma)}{2}}(\rho^b)_x dx\\
=&\frac{2s}{b}(\int_{-\infty}^x \rho^{2s+1-(\alpha+\gamma)}dx)^{\frac{1}{2}}(\int_{-\infty}^x (\rho^b)_x^2 dx)^{\frac{1}{2}}
\leq C\|\rho\|_{L^{\gamma}(R)}\|(\rho^b)_x\|_{L^2(R)}.
\end{aligned}
\end{equation}
Consequently,
\begin{equation}\label{2.2.92'}
\begin{aligned}
\int_{0}^t\sup_{x\in R}\rho^{4s}dt\leq C\sup_{x\in R}\|\rho\|_{L^\gamma(R)}^2
\int_{0}^t\|(\rho^b)_x\|_{L^2(R)}^2 dt
\leq C.
\end{aligned}
\end{equation}
Moreover, applying($\ref{2.2.86'}$), one has
\begin{equation}\label{2.2.93'}
\begin{aligned}
\int_{0}^t\int_{R}(\rho^s)^{4+2l}dxdt
\leq \int_{0}^t(\sup_{x\in R}\rho^{4s}\int_{R}\rho^{2sl}dx)dt
\leq \sup_{t} \int_{R}\rho^{2sl}dx\int_{0}^t\sup_{x\in R}\rho^{4s}dt
\leq C,
\end{aligned}
\end{equation}
where $l\geq 1$ is any real number.

Denote $f(t)=\int_{R}(\rho^s)^{4+2l}dx$. Then, from
($\ref{2.2.86'}$) and ($\ref{2.2.93'}$), one has $f(t)\in L^1(0,
\infty)\cap L^\infty(0, \infty)$. The left is same as in Case 1
($\bar\rho>0$).

The proof of the theorem is finished.
\end{proof}

The proof of Theorem \ref{thm2.2.3} is  completely same as in
\cite{LLX2008,JX2008,JWX2011}and we omit it here.


\begin{thebibliography}{99}


\bibitem{BD2003}
D. Bresch and B. Desjardins, Existence of global weak solutions for
a 2D viscous shallow water equations and convergence to the
quasi-geostrophic model, {\it Comm. Math. Phys.}, {\bf 238}(1-2)
(2003) 211-223.
\bibitem{BD2004}
D. Bresch and B. Desjardins, Quelques modeles diffusifs
capillaires de type Korteweg,
{\it C. R. Acad. Sci. }, Paris, section
mecanique {\bf 332} (11)(2004) 881-886.
\bibitem{BD2006}
D. Bresch and B. Desjardins, On the construction of approximate
solutions for the 2D viscous shallow water model and for
compressible Navier-Stokes models,
{\it J. Math. Pures Appl.} {\bf 86} (2006),
362-368.
\bibitem{BDL2003}
D. Bresch, B. Desjardins, Chi-Kun Lin, On some compressible fluid
models: Korteweg, lubrication, and shallow water systems, {\it
Comm. Partial Differential Equations} {\bf 28}(3-4)( 2003),
843-868.
\bibitem{BDG2007}
D. Bresch, B. Desjardins, D. Gerard-Varet, On compressible
Navier-Stokes equations with density dependent viscosities in bounded domains,
{\it J. Math. Pures Appl.} {\bf 87} 2  (2007)  227-235.
\bibitem{Danchin2000}
R. Danchin, Global existence in critical spaces for compressible
Navier-Stokes equations, {\it Invent. Math.} {\bf 141}(2000),
579-614.
\bibitem{DJ2010}
C. S. Dou, Q. S. Jiu, A remark on free boundary problem of
1-D compressible Navier-Stokes equations
with density-dependent viscosity, {\it Math. Meth. Appl. Sci.} {\bf 33}(2010),
 103¨C116.
\bibitem{FNP2001}
E. Feireisl, A. Novotn\'y and H. Petzeltov\'a,
On the existence of globally defined weak solutions to the
Navier-Stokes equations of isentropic compressible fluids, {\it J.
Math. Fluid Mech.} {\bf 3} (2001) 358-392.
\bibitem{GP2001}
J.F. Gerbeau, B. Perthame. Derivation of Viscous Saint-Venant System for Laminar
Shallow Water; Numerical Validation, Discrete and Continuous Dynamical Systems,
{\it Ser. B, Vol. 1, Num.} {\bf 1} (2001) 89-102.
\bibitem{GJX2008}
Z. Guo, Q. Jiu, Z. Xin, Spherically symmetric isentropic compressible flows with density-dependent
viscosity coefficients, {\it SIAM J. Math. Anal.} {\bf 39} 5 (2008)  1402-1427.
\bibitem{Hoff1987}
D. Hoff, Global existence of 1D compressible isentropic Navier-Stokes equations with large initial data,
{\it Trans. Amer. Math. Soc.} {\bf 303(1)}
(1987) 169-181.
\bibitem{Hoff1995}
D. Hoff, Strong convergence to global solutions for multidimensional flows of compressible viscous fluids with
polytropic equations of state and discontinuous initial data,
{\it Arch. Rat. Mech. Anal.} {\bf 132}
(1995) 1-14.
\bibitem{Hoff1998}
D. Hoff, The zero-Mach limit of compressible flows,
{\it Comm. Math. Phys.} {\bf 192} (1998) 543-554.
\bibitem{HS1991} D. Hoff, D. Serre, The failure of continuous
dependence on initial data for the Navier-Stokes equations of
compressible flow, {\it SIAM J. Appl. Math.} {\bf 51}(1991)
887-898.
\bibitem{HS2001}
D. Hoff, J. Smoller, Non-formation of vacuum states for compressible Navier-Stokes equations.
{\it Comm. Math. Phys.} {bf 216} (2001) no. 2, 255-276.
\bibitem{Jiang1998}
S. Jiang, Global smooth solutions of the equations of a viscous,
heat-conducting one-dimensional gas with density-dependent
viscosity, {\it Math. Nachr.} {\bf 190}(1998) 169-183.
\bibitem{JXZ2005} S. Jiang, Z. P. Xin and P. Zhang, Global weak
solutions to 1D compressible isentropy Navier-Stokes with
density-dependent viscosity,  {\it Methods and Applications of
Analysis} {\bf 12 (3)}(2005) 239-252.
\bibitem{JX2008}
Q. Jiu, Z. P. Xin, The Cauchy problem for 1D compressible flows with density-dependent viscosity coefficients,
{\it Kinet. Relat. Models} {\bf 1} (2) (2008) 313-330.
\bibitem{JWX2011}
Q. Jiu, Y. Wang, Z.P. Xin, Stability of Rarefaction Waves to the 1D Compressible Navier-Stokes Equations with Density-dependent Viscosity,
{\it Comm. Partial Differential Equations} {\bf 36} (2011) 602-634.
\bibitem{JWX2011-2}
Q. Jiu, Y. Wang, Z.P. Xin, Global well-posedness of 2D compressible
Navier-Stokes equations with large data and vacuum,
arXiv:1202.1382v1.
\bibitem{Kanel1968}
J. I. Kanel, A model system of equations for the one-dimensional
motion of a gas. (Russian) {\it Differencial'nye Uravnenija} {\bf 4}
(1968) 721-734.
\bibitem{KS1977}
A. V. Kazhikhov, V. V. Shelukhin, Unique global solution with
respect to time of initial-boundary value problems for
one-dimensional equations of a viscous gas, {\it J. Appl. Math.
Mech. }{\bf 41}(1977) 273-282; translated from {\it Prikl. Mat. Meh.
} {\bf 41 }(1977) 282-291.
\bibitem{LLX2008}
H. L., Li, J. Li, Z. P. Xin, Vanishing of Vacuum States and Blow-up
Phenomena of the Compressible Navier-Stokes Equations, {Comm. Math. Phys.} {\bf 281}(2) (2008) 401-444.
\bibitem{Lions1998}
P. L. Lions, {\it Mathematical Topics in Fluid Dynamics 2,
Compressible Models}, Oxford Science Publication, Oxford, 1998.
\bibitem{LXY1998}
T. P. Liu, Z. P. Xin, T. Yang, Vacuum states of compressible flow,
{\it Discrete Conti
b33
nuous Dynam. systems}{\bf 4}(1)(1998) 1-32.
\bibitem{MN1980}
A. Matsumura, T. Nishida, The initial  value problem for the
equations of motion of viscous and  heat-conductive gases, {\it J.
Math. Kyoto Univ.} {\bf 20}(1) (1980) 67-104.
\bibitem{MV2007}
A. Mellet and A. Vasseur, On the isentropic compressible
Navier-Stokes equation, {\it Comm. Partial
Differential Equations} {\bf
32}(2007) 431-452.
\bibitem{MV2008}
A. Mellet and A. Vasseur, Existence and uniqueness of global strong solutions for one-dimensional
compressible Navier-Stokes equations,
{\it SIAM J. Math. Anal.} {\bf
39}(2008) No. 4, 1344-1365.
\bibitem{Nash1962}
J. Nash, Le probleme de Cauchy pour Les equations differentielles d'un fluids general,
{\it Bulletinde la S. M. F.} {\bf 90} (1962) 487-497.
\bibitem{OMM2002}
M.Okada, S.Matusu-Necasova, T.Makino, Free boundary
problem for the equation of one-dimensional motion of compressible
gas with density-dependent viscosity, {\it Ann. Univ. Ferrara Sez.
VII (N.S.)} {\bf 48} (2002) 1-20.
\bibitem{Serre1986}
D. Serre, Sur l\'equation monodimensionnelle d\'un fluide visqueux, compressible et conducteur de chaleur,
{\it C. R. Acad. Sci. Paris S\'er. I} {\bf 303} (1986) 703-706.
\bibitem{VK1995} V.A. Vaigant and A.V. Kazhikhov, On existence of global
solutions to the two-dimensional Navier-Stokes equations for a
compressible viscous fluid, {\it Siberian J. Math.} {\bf 36} (1995)
1283-1316.
\bibitem{VYZ2003}
S. W. Vong, T. Yang, C. J. Zhu, Compressible Navier-Stokes equations
with degenerate viscosity coefficient and vacuum II, {\it J.
Differential Equations}{\bf 192} (2)(2003) 475-501.
\bibitem{WHG2011}
J. Wei, L. He, Z. Guo, A Remark on the Cauchy Problem of 1D
Compressible Navier-Stokes Equations with Density-dependent
Viscosity Coefficients, Accepted by Acta Mathematicae Applicatae
Sinica, English Series, 2011.
\bibitem{Xin1998}
Z. P. Xin, Blow-up of smooth solution to the compressible
Navier-Stokes equations with compact density, {\it Comm. Pure Appl.
Math.} {\bf 51}(1998) 229-240.
\bibitem{YYZ2001}
T. Yang, Z. A. Yao, C. J. Zhu, Compressible Navier-Stokes equations
with density-dependent viscosity and vacuum , {\it Comm. Partial
Differential Equations} {\bf 26} (5-6)(2001) 965-981.
\bibitem{YZ2002}
T. Yang, H. J. Zhao, A vacuum problem for the one-dimensional
Compressible Navier-Stokes equations with density-dependent
viscosity, {\it J. Differential Equations} {\bf 184} (1)(2002)
163-184.
\bibitem{YZ2002}
T. Yang, C. J. Zhu, Compressible Navier-Stokes equations with
degenerate viscosity coefficient and vacuum ,{\it Comm. Math.
Phys.} {\bf 230} (2)(2002) 329-363.
\bibitem{YYZ2001}
T. Yang, Z. A. Yao, C. J. Zhu, Compressible Navier-Stokes equations
with density-dependent viscosity and vacuum , {\it Comm. Partial
Differential Equations} {\bf 26} (5-6)(2001) 965-981.





\end{thebibliography}
\end{document}